\documentclass{amsart}

\usepackage{amscd,amsmath,amssymb,amsthm,amsfonts,cite,epsfig,graphics}
\usepackage{graphicx}
\usepackage{enumerate}
\usepackage{epstopdf}

\usepackage{mathtools}
\usepackage{caption}
\usepackage{subcaption}
\usepackage{xcolor}

\usepackage{import}
\usepackage{xifthen}
\usepackage{pdfpages}
\usepackage{transparent}
\usepackage{verbatim}

\usepackage{amssymb}
\usepackage{url}
\usepackage{amscd}
\usepackage{texdraw}
\usepackage{tikz-cd}
\usepackage{tikz}

\theoremstyle{plain}
\newtheorem*{starthma}{Theorem A}
\newtheorem*{starthmb}{Theorem B}
\newtheorem*{starthmc}{Theorem C}
\newtheorem*{starthmd}{Theorem D}
\newtheorem*{starthme}{Theorem E}
\newtheorem*{starthmf}{Theorem F}
\newtheorem*{starthmg}{Theorem G}

\newtheorem{theorem}{Theorem}
\newtheorem{corollary}{Corollary}
\newtheorem{lemma}{Lemma}
\newtheorem{remark}{Remark}
\newtheorem{proposition}{Proposition}
\newtheorem{definition}{Definition}
\newtheorem{example}{Example}

\newtheorem{problem}{Problem}
\input txdtools

\newcommand{\tarc}{\mbox{\large$\frown$}}

\newcommand{\arc}[1]{\stackrel{\tarc}{#1}}

\newcommand{\Si}{\Sigma}
\newcommand{\s}{\sigma}
\newcommand{\w}{w=\cdots i_{n-1}\cdots i_{1}i_{0}}
\newcommand{\sw}{\sigma(w)=\cdots i_{n-1}\cdots i_{1}}

\newcommand{\CP}{{\mathbb C}}
\newcommand{\CC}{\widehat{\mathbb C}}
\newcommand{\R}{\mathbb R}
\newcommand{\T}{\mathbb T}
\newcommand{\D}{\mathbb D}
\newcommand{\UH}{\mathbb H}
\newcommand{\HH}{\widehat H}

\begin{document}

\title[Uniqueness and Teichmüller Structure]{Space of Geometric Gibbs Measures: Uniqueness, Ergodicity, and Teichmüller Complex Manifold Structure}

\author{Yunping Jiang}
\address{Department of Mathematics\\
Queens College of the City University of New York\\
Flushing, NY 11367-1597\\
and\\
Department of Mathematics\\
Graduate School of the City University of New York\\
365 Fifth Avenue, New York, NY 10016}
\email{yunping.jiang@qc.cuny.edu}
\thanks{This research is supported partially by the Simons collaboration grants (\#199837, \#523341 and \#942077) for mathematicians and PSC-CUNY awards.}


\begin{abstract}
Building on our previous work, in which geometric Gibbs measures were introduced to complete the space of classical Gibbs measures under the Teichmüller metric, we establish the uniqueness of geometric Gibbs measures in this paper. This result completes the construction of a Teichmüller complex manifold structure on the space of all geometric Gibbs measures. Using the welding technique, we uncover a connection between geometric Gibbs measures and harmonic measures on dynamical quasicircles. This connection allows us to establish asymptotic conformal rigidity for these dynamical quasicircles. Finally, we establish a foundation for the systematic study of the ergodic properties of geometric Gibbs measures.
\end{abstract}

\maketitle

\section{\bf Introduction}~\label{intro}

Teichm\"uller theory is one of the central developments in modern mathematics for understanding the geometry of Riemann surfaces through quasiconformal mappings. Its fundamental object is the \emph{universal Teichm\"uller space}.
The main purpose of this paper is to develop an analogy of Teichm\"uller theory for fractal sets arising in dynamical systems, building on the framework established in our previous work~\cite{JGibbs0,JGibbsIntro,JGibbs}.

We give a brief review of the universal Teichm\"uller space (refer to~\cite{A,AB,AhlforsBook,GL,N,IT,P}). A Beltrami coefficient $b(z)$ on the unit disk $\D$ is a complex-valued measurable function satisfying $\|b\|_\infty<1$. An important problem in the development of the universal Teichm\"uller space used to be solving the Beltrami equation
\[ 
\frac{\partial w}{\partial\overline z} = b(z)\frac{\partial w}{\partial z}.
\] 
Gauss gave the first affirmative answer in the analytic setting through the \emph{isothermal coordinate theorem} (1827), proving that every analytic surface admits local isothermal coordinates. Korn–Lichtenstein later extended the proof to lower regularities. Morrey extended further to substantially weaker regularity by developing powerful methods from elliptic partial differential equations. Building on Vekua's pioneering work, Bojarski established the existence of quasiconformal homeomorphism solutions to the Beltrami equation with measurable Beltrami coefficients. Ahlfors and Bers further developed the theory, establishing, in particular, the holomorphic dependence of normalized solutions on the Beltrami coefficient. These developments culminated in the celebrated \emph{measurable Riemann mapping theorem}, which asserts the existence, uniqueness, and holomorphic dependence, after suitable normalization, of a quasiconformal homeomorphism solving the Beltrami equation. Two Beltrami coefficients are said to be equivalent if their normalized quasiconformal homeomorphism solutions induce the same homeomorphism on the boundary circle of the unit disk. The space of equivalence classes equipped with the Teichm\"uller metric forms the universal Teichm\"uller space. It is an infinite-dimensional complex manifold. 

Motivated by this classical theory, we have started to develop an analogy of Teichm\"uller theory for fractal sets arising in dynamical systems. Consider the symbolic space of sequences of $ 0$'s, $\cdots$, and $d-1$'s with the shift dynamics. In~\cite{JGibbs0,JGibbsIntro,JGibbs}, we introduced the notion of \emph{geometric Gibbs measures} as a metric on the space such that the derivative of the shift dynamics with respect to this metric is a given geometric continuous potential. The main purpose of this paper is to prove that the space of geometric Gibbs measures is the completion of the space of classical Gibbs measures with respect to the Teichm\"uller metric, giving a natural analogy to the universal Teichm\"uller space in symbolic dynamics. We proved the existence in~\cite{JGibbs}. This paper establishes the uniqueness through the following theorem.

\medskip
\begin{starthma}\label{A}
For every geometric continuous $g$-function, there is a unique geometric Gibbs measure associated with it.
\end{starthma}

This uniqueness lets us complete the construction of the Teichm\"uller complex manifold structure on the space of geometric Gibbs measures as proposed in~\cite{JGibbs}. 

\medskip
\begin{starthmb}\label{B}
The space of geometric Gibbs measures is an infinite-dimensional complete complex manifold, which completes the space of classical Gibbs measures under the Teichm\"uller metric. 
\end{starthmb}

We describe the tangent space of this complex manifold at the basepoint in Theorem~C. In Theorem~D, we show that the Hilbert transformation gives an almost complex structure on the tangent space 

Using conformal welding, we weld the uniformly symmetric circle endomorphism associated with a geometric Gibbs measure to the polynomial map $z\mapsto z^{d}$, thereby introducing the notions of a dynamical quasicircle (a fractal set in general) and a dynamical quasidisk. The harmonic measure on a dynamical quasicircle, viewed from the center $0$, is a geometric Gibbs measure on the quasicircle.

\medskip
\begin{starthme}\label{E}
The harmonic measure on a dynamical quasicircle centered at $0$ is a geometric Gibbs measure.
\end{starthme}

Applying Theorem~A, we obtain the following asymptotic rigidity result.

\medskip
\begin{starthmf}\label{G}
If two dynamical quasicircles are asymptotically conformally equivalent, then they coincide up to conjugations by M\"obius transformations.
\end{starthmf}

Unlike classical Gibbs measures, which are always ergodic, geometric Gibbs measures need not be ergodic. We therefore initiate a systematic study of the ergodic theory of geometric Gibbs measures. To this end, we introduce the notion of a uniformly absolutely continuous circle endomorphism and prove the following result.

\medskip
\begin{starthmg}\label{H}
The unique geometric Gibbs measure associated with a geometric continuous $g$-function is ergodic whenever its induced circle endomorphism is uniformly absolutely continuous.
\end{starthmg}

The uniform absolute continuity may have a connection with the Weil--Petersson geometry. For this purpose, we define a uniform Weil--Petersson property for the induced circle endomorphism for a geometric Gibbs measure (see Remark~\ref{uwp}). 

The paper is organized as follows. In Section~\ref{geomea}, we review the definition of a geometric Gibbs measure associated with a geometric continuous $g$-function, together with several results established in~\cite{JGibbs}. We state these results without proofs and refer the reader to~\cite{JGibbs} for details. Section~\ref{uniq} proves the uniqueness theorem (Theorem~A). In Section~\ref{ccus}, we establish a complex characterization of uniformly symmetric circle endomorphisms, which, together with Theorem~A, enables us to complete the construction of the Teichm\"uller complex manifold structure in Section~\ref{tcms} (Theorem~B and Corollary~\ref{completion2}). In the same section, we describe the tangent space of this complex manifold at the basepoint (Theorem~C) and the almost complex structure (Theorem~D). In Section~\ref{weld}, we introduce dynamical quasicircles and dynamical quasidisks via conformal welding with the polynomial map $z\mapsto z^{d}$, prove that the harmonic measure of a dynamical quasidisk is a geometric Gibbs measure (Theorem~E), and establish the asymptotic conformal rigidity theorem (Theorem~F). Finally, in Section~\ref{erg}, we begin studying the ergodicity of a geometric Gibbs measure by introducing uniformly absolutely continuous circle endomorphisms and proving Theorem~G. We conclude the paper with a remark (Remark~\ref{uwp}), where we discuss a possible connection between the uniform absolute continuity and the uniform Weil--Petersson geometry.

\section{\bf Geometric Gibbs Measures}~\label{geomea}

In~\cite{JGibbs0,JGibbsIntro,JGibbs}, we introduced geometric Gibbs measures within the Teichm\"uller metric framework. In this section, we review this definition and highlight key results from~\cite{JGibbs}.

Let $d\geq 2$ be a fixed integer. Consider the symbolic space, 
$$
\Si=\prod_{\infty}^{0}\{0, \ldots, d-1\}=
$$
$$
\Big\{ \w \; |\; i_{n-1}\in \{0,\cdots, d-1\}, n\geq 1\Big\}
$$
with the product topology. This space is compact and metrizable. The right shift
$$
\s : \w \to \sw
$$
defines a symbolic dynamical system on $\Si$. 

For a point $\w\in \Si$, the (right) $n$-cylinder $[w]_{n}$ containing $w$ is the subset
$$
[w]_{n}=\{ w'= \cdots i_{n+m}'\cdots i_{n}'i_{n-1}\cdots i_{0}\in \Si\;|\;  i_{n+m}' \in \{0, \cdots, d-1\},\; m\geq 0\}.
$$ 
The set ${\mathcal B}_{0}$ of all cylinders forms a topological basis for $\Si$.

Let ${\mathcal C}={\mathcal C} (\Si)$ denote the space of all continuous real-valued functions on $\Si$.  
A function $\phi\in \Si$ is called H\"older continuous if there exist constants 
$C>0$ and $0<\tau < 1$ such that for all $w\in \Si$, $n\geq 1$, and $w' \in [w]_{n}$,
$$ 
|\phi (w) -\phi (w')| \leq C \tau^{n}.
$$
Let ${\mathcal C}^{H}={\mathcal C}^{H} (\Si)$ be the space of all H\"older continuous functions. 
We refer to a positive function in $\mathcal{C}$ as a {\em continuous potential}, and a positive function in 
$\mathcal{C}^{H}$ as a {\em smooth potential}.
 
 Let ${\mathcal B}$ denote the standard Borel $\sigma$-algebra generated by ${\mathcal B}_{0}$.
Let ${\mathcal M}$ be the space of all finite Borel measures on $\Si$. By the Riesz representation theorem, ${\mathcal M}=\mathcal{C}^{*}$, the dual space of $\mathcal{C}$. 
A measure $\mu\in \mathcal{M}$ is called a probability measure if $\mu (\Si)=1$.  Let ${\mathcal M}_{1}$ be the set of all probability measures in $\mathcal{M}$; this set is convex and compact under the weak $*$-topology. 

A measure $\mu\in {\mathcal M}_{1}$ is said to be $\sigma$-invariant if
$$
\mu (\sigma^{-1}(A)) =\mu (A), \quad  \forall A\in {\mathcal B}.
$$ 
Let ${\mathcal M}_{1}(\sigma)$ be the set of all $\sigma$-invariant measures in ${\mathcal M}_{1}$. This set is non-empty, convex, and compact under the weak $*$-topology. Define $\mu_{0}$ as the measure such that $\mu_{0} ([w]_{n})=1/d^{n}$ for all $w\in \Sigma$. Then $\mu_{0}\in {\mathcal M}_{1}(\sigma)$ and is a Bernoulli measure. 

The classical Gibbs theory developed by Sinai and Ruelle (see~\cite{Sinai1,Sinai2, Ruelle1, Ruelle2, Ruelle3}) guarantees that, for a given smooth potential $\psi$,  there is a number $P=P(\log \psi)$, known as the pressure, and 
a unique measure $\mu= \mu_{\psi}\in {\mathcal M}_{1}(\sigma)$, called the Gibbs measure associated with $\psi$, and a constant $C>0$ such that
$$
C^{-1} \leq \frac{\mu ([w]_{n})}{\exp (-Pn + \sum_{k=0}^{n-1} \log
\psi (\sigma^{k} (w)))}\leq C, \quad \forall \;[w]_{n}\in {\mathcal B}_{0}.
$$
Note that $e^{P}$ is the unique maximal eigenvalue of Ruelle's Perron-Frobenius operator associated with the potential $\psi$ (refer to~\cite{JR}).  

Two smooth potentials $\psi_{1}$ and $\psi_{2}$ are said to be cohomologous, denoted by $\psi_{1} \sim_{ch} \psi_{2}$,  
if there is a continuous function $u$ such that
$$
\log \psi_{2} = \log \psi_{1} + u\circ \sigma - u.
$$
If $\psi_{1}\sim_{ch} \psi_{2}$, then the Gibbs measures associated with them are the same, i.e., $\mu_{\psi_{1}}=\mu_{\psi_{2}}$. Converesly, if $\mu_{\psi_{1}}=\mu_{\psi_{2}}$, then $\psi_{1}\sim_{ch} \psi_{2}$. Given $\psi\in \mathcal{C}^{H}$, let 
$$
[\psi] =\{ \widetilde{\psi} \sim_{ch} \psi\;|\; \widetilde{\psi}\in \mathcal{C}^{H}\}
$$
be the cohomologous equivalence class of $\psi$. 
There is a bijective map between the space of all cohomologous equivalence classes and 
the space of all Gibbs measures. From Ruelle's Perron-Frobenius theorem (refer to~\cite{JR}), for each $\psi\in \mathcal{C}^{H}$, there is a unique function $g\in \mathcal{C}^{H}$ with $e^{P} g\in [\psi]$ satisfying  
\begin{equation}~\label{sum}
\sum_{i=0}^{d-1} g(wi) =1, \quad \forall w\in \Si.
\end{equation}
Traditionally, such a potential is called a (smooth) $g$-function (see~\cite{K,Le} and~\cite[Theorem~2.2]{JGibbs}).

The constant function $g_{0}\equiv 1/d$ is a $g$-function, and its corresponding Gibbs measure is $\mu_{0}$. However, a typical (smooth) $g$-function is neither constant nor even locally constant. Moreover, the space of all smooth $g$-functions is in bijective correspondence with the space of all Gibbs measures.

In~\cite{JGibbs0,JGibbsIntro,JGibbs}, we concentrated on the space of (smooth) $g$-functions
$$
\mathcal{HG} =\{ g_{f} \:|\; f\in {\mathcal E}^{1+H}\},
$$
where ${\mathcal E}^{1+H}$ is the space of all $C^{1+\alpha}$ expanding circle endomorphisms $f$ of degree $d$ for some 
$0<\alpha\leq 1$. When $d=2$, we know that $g\in \mathcal{HG}$ if and only if it satisfies a compatibility condition (refer to~\cite{CJQ}). Let us briefly describe $\mathcal{HG}$, with more detail in~\cite{JGibbs}. 

Let $\CP$ denote the complex plane and let $\D=\{ z\in \CP\;|\; |z|<1\}$ be the open unit disk in $\CP$. Suppose  
$$
\T=\partial \D= \{ z\in \CP\; |\; |z|=1\}
$$
represents the unit circle. The real line $\R$ serves as the universal cover of $\T$. The universal covering map is given by  
$$
\pi (x) = e^{2\pi i x}: \R \to \T
$$
with $\pi (0)=1$. 

Consider an orientation-preserving covering map of degree $d$
$$
f: \T\to \T.
$$
Every such map $f$ can be lifted to an
orientation-preserving homeomorphism
$$
F: \R\to \R, \quad F(x+1)=F(x)+d, \quad \forall
x\in {\mathbb R}.
$$
Throughout this paper, we assume $F(0)=0$, which implies $f(1)=1$. 
Under this assumption, $f$ and $F$ correspond one-to-one. We refer to $F$ as the lift of $f$, and either $f$ or $F$ as a circle endomorphism of degree $d$. 

We say that $f$ is $C^{1+\alpha}$ for some $0<\alpha\leq 1$ if its derivative $F'$ exists and is an $\alpha$-H\"older continuous function. We say that $f$ is expanding if there exist constants $0<C\leq 1$ and $\lambda >1$ such that
$$
(F^{n})'(x) \geq C\lambda^{n}, \quad x\in \R, \quad \forall n\geq 1.
$$
We use $\mathcal{E}^{1+H}$ to denote the space of all $C^{1+\alpha}$, $0<\alpha\leq 1$, expanding circle endomorhisms $f$ of degree $d$. 

We also consider an orientation-preserving circle homomorphism,
$$
h: {\mathbb T}\to {\mathbb T}.
$$ 
Let $H: {\mathbb R}\to {\mathbb R}$ be its lift of $h$. Then, we have 
$$
H(x+1)=H(x)+1, \quad \forall\;\; x\in {\mathbb R}.
$$  
Throughout this paper, we will assume that $H(0)=0$, which implies $h(1)=1$.
We refer to $H$ as the lift of $h$ and either $h$ or $H$ as a circle homeomorphism. 

For any interval $I=[a,b]$ in $[0, 1]$, we denote its length by $|I|=b-a$. This notation also applies to the length $|J|$ of the interval $J\subset \T$ with $\pi (I)=J$. We use $|\cdot |$ to denote the Lebesgue probability measure on $\T$. The measure $|\cdot |$ on any interval $J\subset \T$ is $|J|=|I|$. 

Consider another symbolic space defined as 
$$
\Si^{+}=\prod_{0}^{\infty} \{ 0, 1, \cdots, d-1\}
=
$$
$$
\{ v=j_{0}j_{1}\cdots j_{n-1} \cdots \; |\; j_{n-1}\in
\{ 0, 1, \cdots, d-1\}, \; n\geq 1\}
$$
with the product topology. The left shift 
$$
\s^{+} (v) =  j_{1}\cdots j_{n-1}j_{n} \cdots
$$
defines a symbolic dynamical system on $\Si^{+}$.

For a point $v=j_{0}j_{1}\cdots j_{n-1}\cdots \in \Si^{+}$, the (left) $n$-cylinder containing $v$ is the subset 
$$
[v]_{n}^{+} =\{ v' =j_{0}j_{1}\cdots j_{n-1}j_{n}'j_{n+1}' \cdots
\;|\; j_{n+k}' \in \{ 0, 1, \cdots, d-1\}, k=0, 1, \cdots\}
$$
The set $\mathcal{B}_{0}^{+}$ of all cylinders forms a topological basis for $\Sigma^{+}$. Let $\mathcal{B}^{+}$ be the standard Borel $\s$-field generated by $\mathcal{B}_{0}^{+}$.

Let ${\mathcal M}^{+}$ be the space of all Borel finite measures on $\Sigma^{+}$, and let ${\mathcal M}^{+}_{1} (\sigma^{+})$ denote the space of all $\sigma^{+}$-invariant probability measures in $\mathcal{M}^{+}$. Specifically, these are the probability measures $\mu^{+}\in {\mathcal M}^{+}$
such that 
$$
\mu^{+} ((\sigma^{+})^{-1} (A)) =\mu^{+} (A), \quad \forall\;\; A\in \mathcal{B}^{+}.
$$  

The space $\Sigma^{+}$ can be thought of as the symbolic representation of the unit circle $\T$. Specifically, any point $z=e^{2\pi i x}\in \T$ can be represented by a point $v=j_{0}j_{1}\cdots j_{n-1} \cdots\in \Sigma^{+}$ for 
\begin{equation}~\label{rep}
x=x(v) =\sum_{k=0}^{\infty} \frac{j_{k}}{d^{k+1}}.  
\end{equation}
We define the length on $\Si^{+}$ as $|v-v'| =|x(v)-x(v')|=|z-z'|$. 
The Lebesgue measure $|\cdot |$ on $\T$ can be lifted to $\Si^{+}$ as a probability measure $\mu_{0}^{+}$ with $\mu_{0}^{+} ([v]_{n}^{+}) =1/d^{n}$. Then $\mu_{0}^{+}\in {\mathcal M}_{1}^{+}(\sigma^{+})$ and is a Bernoulli measure. 

For a circle endomorphism $f$ of degree $d$, 
the preimage $f^{-1}(1)$ divides $\T$ into $d$ closed
intervals $J_{0}$, $J_{1}$, $\cdots$, $J_{d-1}$, ordered in a
counter-clockwise direction. If $J_{0}$ has an endpoint
$1$, then $J_{d-1}$ also has $1$ as an endpoint. We define
$$
\varpi_{0}=\varpi_{0,f}=\{ J_{0}, J_{1}, \cdots, J_{d-1}\}.
$$
This set serves as a Markov partition, satisfying the following conditions:
\begin{itemize}
\item[{\rm i.}] $\T=\cup_{j=0}^{d-1} J_{j}$,
\item[{\rm ii.}] the restriction of $f$ to the interior of $J_{j}$ is
injective for every $0\leq j\leq d-1$,
\item[{\rm iii.}] $f(J_{j})=\T$ for every $0\leq j\leq d-1$.
\end{itemize}

Let $I_{0}$, $I_{1}$, $\cdots$, $I_{d-1}$ be the lifts of $J_{0}$,
$J_{1}$, $\cdots$, $J_{d-1}$ in $[0,1]$. We have the following properties:
\begin{itemize}
\item[{\rm i)}] $[0,1]=\cup_{j=0}^{d-1} I_{j}$,
\item[{\rm ii)}] $F(I_{j})=[j, j+1]$ for every $0\leq j\leq d-1$.
\end{itemize}
We denote the set of these intervals as 
$$
\eta_{0}=\eta_{0,f}=\{I_{0}, I_{1}, \cdots, I_{d-1}\}.
$$
Thus $\eta_{0}$ constitues a Markov partition of $[0,1]$.

Let $v_{n}=j_{0}j_{1} \cdots j_{n-1}$ be a
word of $0's$, $1's$, $\cdots$, and $(d-1)'s$ of length $n$. Then we have a unique interval $J=J_{v_{n}}$ such that  
\begin{enumerate}
\item $f^{k}(J) \subset J_{j_{k}}$ for $0\leq k\leq n-2$,
\item $f^{n-1} :J \to J_{j_{n-1}}$ is a homomorphism. 
\end{enumerate}
The partition 
$$
\varpi_{n}=\varpi_{n,f} 
$$
$$
=\{ J_{v_{n}}\;|\; v_{n}
=j_{0}j_{1}\cdots j_{n-1},\;
j_{k}\in\{ 0, 1, \cdots, d-1\},\; k=0, 1,\cdots, d-1\}.
$$
is considered as the pull-back partition of $\varpi_{0}$ by $f^{n}$. It contains $d^{n}$ closed intervals and also serves as 
a Markov partition of $\T$. 
Let $\eta_{n}$ be the corresponding lift partition of $\varpi_{n}$
in $[0,1]$ with the same labeling. Denote
$$
\eta_{n} =\eta_{n,f} 
$$
$$
=\{ I_{v_{n}}\;|\; v_{n}=j_{0}j_{1}\cdots j_{n-1},\;
j_{k}\in\{ 0, 1, \cdots, d-1\},\; k=0, 1,\cdots, d-1\}.
$$
We refer to the sequences of partitions $\{\varpi_{n}\}_{n=0}^{\infty}$ or $\{\eta_{n}\}_{n=0}^{\infty}$
as
the sequence of nested Markov partitions generated by $f$. 

For a point $v=j_{0}\cdots j_{n-1}\cdots \in \Sigma^{+}$, we define the truncated word as 
$v_{n}:=j_{0}\cdots j_{n-1}$, which allows us to establish a nested sequence of intervals:
$$
\cdots \subset J_{v_{n}} \subset J_{v_{n-1}} \subset \cdots
J_{v_{1}}\subset \T.
$$
Since each $J_{v_{n}}$ is compact,
$$
J_{v}=\cap_{n=1}^{\infty} J_{v_{n}} \neq \emptyset.
$$
When every $J_{v}=\{ x_{v}\}$ contains only one point (which we assume throughout this paper), we define 
the projection $\pi_{f}$ from $\Sigma^{+}$ onto $\T$ as
$$
\pi_{f} (v) =x_{v}.
$$
The projection $\pi_{f}$ is one-to-one except for a countable set of points that are mapped to endpoints of all intervals $J_{v_{n}}$, on which $\pi_{f}$ is two-to-one. From our construction, one can verify that
\begin{equation}~\label{sconj}
\pi_{f}\circ \sigma^{+}(v) = f \circ \pi_{f} (v), \quad v\in \Sigma^{+}.
\end{equation}
That indicates that $f$ is always semi-conjugate to $\sigma^{+}$. 
Thus, we refer to $\sigma^{+}: \Si^{+}\to \Si^{+}$ as the topological representation for all $f\in {\mathcal E}^{1+H}$. 

We set up a relationship between $\Si$ and $\Si^{+}$ through $f$. Consider $w=\cdots i_{n-1}\cdots i_{1}i_{0}\in \Sigma$. Let 
$$
w_{n} =i_{n-1}\cdots i_{1}i_{0}\quad \hbox{and}\quad w_{n-1}'=\sigma(w_{n}) = i_{n-1}\cdots i_{1}.
$$
Define $j_{0}=i_{n-1}$, $\cdots$, $j_{n-2}=i_{1}$, $j_{n-1}=i_{0}$. Then we get 
$$
v_{n} =j_{0}\cdots j_{n-2}j_{n-1}\quad \hbox{and}\quad v_{n-1}' = j_{0} \cdots j_{n-2}.
$$
Without cause of confusion, we use $\sigma(w_{n})$ to denote $v_{n-1}'$. Note that $\sigma^{+}(v_{n}) =j_{1}\cdots j_{n-1}$. Then   
$$
J_{w_{n}}=J_{v_{n}} \in \varpi_{n}\quad \hbox{and}\quad J_{\sigma(w_{n})}=J_{v_{n-1}'}\in \varpi_{n-1}
$$ 
and 
$$
J_{w_{n}} \subset J_{\sigma(w_{n})}.
$$
With these notations, define
\begin{equation}~\label{precp}
g_{n} (w) =\frac{|J_{w_{n}}|}{|J_{\sigma(w_{n})}|}.
\end{equation} 
It follows that
\begin{equation}~\label{sumprecp}
\sum_{i=1}^{d-1} g_{n}(wi) =\sum_{i=1}^{d-1} \frac{|J_{w_{n-1}i}|}{|J_{w_{n-1}}|}=1
\end{equation} 
We have the following theorem (see~\cite[Theorem 6.4]{JGibbs}). 

\medskip
\begin{theorem}~\label{sg} Suppose $f\in \mathcal{E}^{1+H}$. Then for every $w \in \Sigma$,
\begin{equation}~\label{cp}
g(w) = \lim_{n\to \infty} g_{n} (w)
\end{equation} 
exists, and the convergence is exponential. Furthermore, the function
$$
g(w)=g_{f}(w): \Sigma \to {\mathbb R}^{+} =(0, \infty)
$$
is a smooth $g$-function.
\end{theorem}

\medskip
\begin{remark}
To compare with our later definition of a geometric continuous $g$-function, we call a $g$-function obtained from Theorem~\ref{sg} a geometric smooth $g$-function.
\end{remark}

Then $\mathcal{HG}$ is the space of all $g_{f}$ for $f\in \mathcal{E}^{1+H}$. One of the most exciting results in~\cite{JGibbs0,JGibbsIntro,JGibbs} is the definition of the Teichm\"uller metric on $\mathcal{HG}$. Another work to study the Teichm\"uller compactification of the space of smooth expanding circle endomorphisms is McMullen's paper ~\cite{M}. 

Consider $q_{d}(z)=z^{d}$ as the basepoint of ${\mathcal E}^{1+H}$. The equation (\ref{sconj}) implies that, for any $f\in {\mathcal E}^{1+H}$, we have a circle homeomorphism $h_{f}$ such that
$$
f\circ h_{f}=h_{f}\circ q_{d}.
$$ 
Thus, we consider a pair $(f, h_{f})$. Two pairs $(f_{1}, h_{f_{1}})$ and $(f_{2}, h_{f_{2}})$ are said to be the Teichm\"uller equivalent if $h_{f_{2}}\circ h_{f_{1}}^{-1}$ is a $C^{1}$-diffeomorphism. From~\cite{JGibbs0,JGibbsIntro,JGibbs}, we know that $g_{f_{1}}=g_{f_{2}}$ if and only if $(f_{1}, h_{f_{1}})$ and $(f_{2}, h_{f_{2}})$ are Teichm\"uller equivalent. 
Use $\tau=[(f, h_{f})]$ to denote the equivalence class of $(f, h_{f})$.
Let $\mathcal{TE}^{1+H}$ denote the space of all Teichm\"uller equivalence classes $\tau$. Then we have a bijective map between $\mathcal{HG}$ with the basepoint $g_{0}=1/d$ and $\mathcal{TE}^{1+H}$ with the basepoint $\tau_{0}=[q_{d}]$:
$$
\zeta: \mathcal{TE}^{1+H} \to \mathcal{HG}; \quad \zeta(\tau)=g_{\tau},
$$
where $g_{\tau}=g_{f}$ for any $(f,h_{f})\in \tau$.

As seen in~\cite[Corollary 5.6, Example 3.10]{JGibbs}), for any pair of $\tau_{1}=[(f_{1}, h_{f_{1}})]$ and $\tau_{2}=[(f_{2}, h_{f_{2}})]$ in $\mathcal{TE}^{1+H}$,  the map $h=h_{f_{2}}\circ h_{f_{1}}^{-1}$ is quasisymmetric, meaning that there exists a constant $M\geq 1$ such that 
\begin{equation}~\label{qsi}
M^{-1} \leq \frac{|H(x+t)-H(x)|}{|H(x)-H(x-t)|} \leq M, \quad \forall x\in {\mathbb R}, \;\; t>0.
\end{equation}
We also call $h$ $M$-quasisymmetric. 

The quasisymmetric property allows us to introduce the Teichm\"uller metric $d_{\mathcal{T}}(\tau_{1}, \tau_{2})$ as follows. 

For any neighborhood $U$ of the unit circle $\mathbb{T}$, the map
$h$ has a quasiconformal extension $\widetilde{h}$ on $U$ with quasiconformal dilatation
$$
K (\widetilde{h},U) = \frac{1+\|b_{\widetilde{h}}\|_{\infty}}{1-\|b_{\widetilde{h}}\|_{\infty}}
$$
where $b_{\widetilde{h}} =\widetilde{h}_{\overline{z}}/\widetilde{h}_{z}\in L^{\infty} (U)$ is the Beltrami coefficient of $\widetilde{h}$ on $U$. 
The boundary quasiconformal dilatation is defined as
$$
B_{h}=\inf_{\widetilde{h}, U} K(\widetilde{h}, U),
$$ 
where the infimum runs over all neighborhoods $U$ and all quasiconformal extensions $\widetilde{h}$ on $U$. This quantity $B_{h}$ is notably independent of specific representatives chosen in $\tau_{1}$ and $\tau_{2}$. 
The Teichm\"uller metric on $\mathcal{TE}^{1+H}$ is defined as 
\begin{equation}~\label{metric2}
d_{\mathcal{T}} (\tau_{1},\tau_{2}) =\frac{1}{2} \log B_{h}.
\end{equation}
This further defines the Teichm\"uller metric on $\mathcal{HG}$ as
$$
d_{\mathcal{HG}} (g_{\tau_{1}}, g_{\tau_{2}}) = d_{\mathcal{T}} (\zeta^{-1} (g_{\tau_{1}}), \zeta^{-1}(g_{\tau_{2}})), \quad g_{\tau_{1}}, g_{\tau_{2}}\in \mathcal{HG}.
$$
Let $\mathcal{HGM}$ be the space of all Gibbs measures associated with $g_{\tau} \in \mathcal{HG}$ with the basepoint $\mu_{0}$. We have a bijective map
$$
\eta: \mathcal{HG} \to \mathcal{HGM}; \quad \eta (g_{\tau}) =\mu_{g_{\tau}}.
$$
We define the Teichm\"uller metric on $\mathcal{HGM}$ with the basepoint $\mu_{0}$ as 
$$
d_{\mathcal{HGM}} (\mu_{\tau_{1}}, \mu_{\tau_{2}}) = d_{\mathcal{HG}} (\eta^{-1} (\mu_{\tau_{1}}), \eta^{-1}(\mu_{\tau_{2}})).
$$
The metric space $(\mathcal{TE}^{1+H}, d_{\mathcal{T}})$,  $(\mathcal{HG}, d_{\mathcal{HG}})$, and $(\mathcal{HGM}, d_{\mathcal{SGM}})$ are all incomplete. We studied their completions 
$(\mathcal{TUS}, d_{\mathcal{T}})$ and $(\mathcal{G}, d_{\mathcal{G}})$ in~\cite{JGibbs0,JGibbsIntro,JGibbs}. Here, we briefly review these completions; see~\cite {JGibbs} for more detail. 

\medskip
\begin{definition}~\label{usce}
A circle endomorphism $f$ is called uniformly symmetric if there
is a bounded positive real function $\varepsilon (t)$ with
$\varepsilon (t) \to 0^{+}$ as $t\to 0^{+}$ such that
$$
\frac{1}{1+\varepsilon (t)} \leq
\frac{|F^{-n}(x+t)-F^{-n}(x)|}{|F^{-n}(x)-F^{-n}(x-t)|}\leq
1+\varepsilon (t)
$$
for all $x\in {\mathbb R}$, all $t>0$, and all $n\geq 1$. 
\end{definition}
 
 Let $\mathcal{US}$ be the set of all uniformly symmetric circle endomorphisms of degree $d$.  From~\cite[Example 3.10]{JGibbs}, we have that
 $$
 \mathcal{E}^{1+H} \subset \mathcal{US}.
 $$ 
Consider $q_{d}(z)=z^{d}$ as the basepoint of $\mathcal{US}$ again.
For any $f\in \mathcal{US}$, let $h_{f}$ be the conjugacy from $f$ to $q_{d}$ (see~\cite[Lemma 4.1]{JGibbs}), which stisfies 
$$
f\circ h_{f}=h_{f}\circ q_{d}.
$$ 
From~\cite[Lemma 5.4 and Lemma 5.5]{JGibbs}, we know that $h_{f}$ is quasisymmetric as defined in (\ref{qsi}). We consider the pair $(f, h_{f})$ for every $f\in\mathcal{US}$. We need the concept of symmetry. 

\medskip
\begin{definition}~\label{sch}
A circle homeomorphism $h$ is called symmetric if there exists a bounded positive real function
$\varepsilon (t)$ with $\varepsilon (t) \to
0^{+}$ as $t\to 0^{+}$ such that
$$
\frac{1}{1+\varepsilon (t)} \leq
\frac{|H(x+t)-H(x)|}{|H(x)-H(x-t)|}\leq 1+\varepsilon (t), \quad \forall \; x\in \R, \;t>0.
$$
\end{definition}

\medskip
\begin{remark}~\label{cmp} Suppose ${\mathcal S}$ is the space of all symmetric circle homeomorphism. Then ${\mathcal S}$ is closed under composition, that is, if $h_{1}$ and $h_{2}$ are both in ${\mathcal S}$, then $h_{2}\circ h_{1}$ is also in ${\mathcal S}$ (refer to~\cite{GS} or~\cite[Lemma 3.6]{JGibbs}).
\end{remark}

Two pairs $(f_{1}, h_{f_{1}})$ and $(f_{2}, h_{f_{2}})$ for $f_{1}, f_{2}\in \mathcal{US}$ are said to be Teichm\"uller equivalent if $h=h_{f_{2}}\circ h_{f_{1}}^{-1}$ is symmetric. Let $\tau=[(f, h_{f})]$ be the Teichm\"uller equivalence class of the pair $(f, h_{f})$. The Teichm\"uller space $\mathcal{TUS}$ is the space of all Teichm\"uller equivalence classes $\tau$ of pairs $(f, h_{f})$ for all $f\in \mathcal{US}$ with the basepoint $\tau_{0}=[q_{d}]$. 
The following theorem is proved in~\cite[Theorem 6.4]{JGibbs}. 

\medskip
\begin{theorem}~\label{sg2} 
Suppose $f\in \mathcal{US}$. Then for every $w \in \Sigma$, the limit
\begin{equation}~\label{cp1}
g(w) = \lim_{n\to \infty} g_{n} (w)
\end{equation} 
exists, and the convergence is uniform. The resulting function
$$
g(w)=g_{f}(w): \Sigma \to {\mathbb R}^{+}
$$
is a continuous $g$-function.
\end{theorem}

\medskip
\begin{remark}
We call a $g$-function obtained from Theorem~\ref{sg2} a geometric continuous $g$-function. 
\end{remark}

\medskip
\begin{remark}
The main difference between Theorem~\ref{sg} and Theorem~\ref{sg2} is the rate of convergence: in Theorem~\ref{sg}, convergence is exponential, whereas in Theorem~\ref{sg2}, it is only uniform. This distinction gives rise to the main difficulty in establishing the uniqueness of a geometric Gibbs measure associated with a continuous geometric $g$-function.
\end{remark}

Let $\mathcal{G}$ represent the space of all $g_{f}$ for $f\in \mathcal{US}$. When $d=2$, $g\in \mathcal{G}$ if and only if it satisfies a compatibility condition (refer to~\cite{CGJ}). 
According to~\cite[Theorem 6.7]{JGibbs}, each $g_{f}\in \mathcal{G}$ serves as a complete symmetric invariant. Specifically, $g_{f_{1}} =g_{f_{2}}$ if and only if $(f_{1}, h_{f_{1}})$ and $(f_{2}, h_{f_{2}})$ are Teichm\"uller equivalent. We get a bijective mapping, which we still denote as $\zeta$,
$$
\zeta: \mathcal{TUS}  \to \mathcal{G} , \quad \zeta(\tau)=g_{\tau},
$$
where $g_{\tau}=g_{f}$ for any $(f, h_{f}) \in \tau$. Furthermore,  
from~\cite[Theorem 7.1]{JGibbs}, we know that
$$
\mathcal{TE}^{1+H} \subset \mathcal{TUS}.
$$
Now, for any pair of $\tau_{1}=[(f_{1}, h_{f_{1}})]$ and $\tau_{2}=[(f_{2}, h_{f_{2}})]$ in $\mathcal{TUS}$, the map 
$h=h_{f_{2}}\circ h_{f_{1}}^{-1}$ is quasisymmetric as shown in~\cite[Corollary 5.6]{JGibbs}. We define the Teichm\"uller distance between $\tau_{1}$ and $\tau_{2}$ using the boundary quasiconformal dilatation of $h$:
for any neighborhood $U$ of the unit circle $\mathbb{T}$, 
$h$ has a quasiconformal extension $\widetilde{h}$ on $U$ with quasiconformal dilatation
$$
K(\widetilde{h}, U) = \frac{1+\|b_{\widetilde{h},U}\|_{\infty}}{1-\|b_{\widetilde{h},U}\|_{\infty}}
$$
where $b_{\widetilde{h},U} =\widetilde{h}_{\overline{z}}/\widetilde{h}_{z}\in L^{\infty} (U)$ is the Beltrami coefficient of $\widetilde{h}$ on $U$. 
The boundary quasiconformal dilatation of $h$ is 
$$
B_{h}=\inf_{\widetilde{h}, U} K (\widetilde{h}, U)
$$ 
where the infimum runs over all neighborhoods $U$ and all quasiconformal extensions $\widetilde{h}$ on $U$. This quantity is notably independent of specific representatives chosen in $\tau_{1}$ and $\tau_{2}$. 
The Teichm\"uller metric on $\mathcal{TUS}$ is defined as 
\begin{equation}~\label{metric3}
d_{\mathcal{T}} (\tau_{1},\tau_{2}) =\frac{1}{2} \log B_{h}.
\end{equation}
This induces the Teichm\"uller metric on $\mathcal{G}$ with the basepoint $g_{0}=1/d$,
$$
d_{\mathcal{G}} (g_{\tau_{1}}, g_{\tau_{2}}) = d (\zeta^{-1} (g_{\tau_{1}}), \zeta^{-1}(g_{\tau_{2}})), \quad g_{\tau_{1}}, g_{\tau_{2}}\in \mathcal{G}.
$$
Furthermore, we prove the following theorem in~\cite[Theorem 8.1]{JGibbs}.

\medskip
\begin{theorem}~\label{completion}
The Teichm\"uller space $(\mathcal{TUS}, d_{\mathcal{T}})$ (resp. $(\mathcal{G}, d_{\mathcal G})$) 
is complete and serves as the completion of $(\mathcal{TE}^{1+H}, d_{\mathcal{T}})$ (resp. $(\mathcal{HG}, d_{\mathcal{HG}})$). That is, 
$$
(\mathcal{TUS}, d_{\mathcal{T}})=\overline{(\mathcal{TE}^{1+H}, d_{\mathcal{T}})}
$$ 
$$
(resp. \quad (\mathcal{G}, d_{\mathcal G})=\overline{(\mathcal{HG}, d_{\mathcal{HG}})}).
$$
\end{theorem}

Since $\mathcal{G}$ is a subspace of the space of all continuous functions on $\Si$ equipped with the maximum norm
$$
\| g \| = \max_{w\in \Si} |g(w)|, \quad \forall g\in \mathcal{G}.
$$
This norm defines the maximum metric on $\mathcal{G}$:
$$
d_{max} (g_{1}, g_{2}) = \| g_{1}-g_{2}\|, \quad g_{1}, \; g_{2} \in \mathcal{G}.
$$
The space ${\mathcal G}$, under $d_{max}$, is incomplete. Thus, this metric is not equivalent 
to the Teichmüller metric $d_{\mathcal{G}}$. However, as shown in~\cite[Theorem 9.1]{JGibbs}, both metrics induce the same topology on ${\mathcal G}$. (Consequently, $d_{\mathcal{G}}$ ais a global hyperbolic metric; in contrast, $d_{max}$ is a local Euclidean metric, like the hyperbolic metric and the Euclidean metric on the open unit disk in the plane.) 

\medskip
\begin{theorem}~\label{unifcont}
The identity map
$$
id_{TM}: (\mathcal{G}, d_{\mathcal{G}}) \to (\mathcal{G}, d_{\max})
$$
is uniformly continuous, whereas the identity map
$$
id_{MT}: (\mathcal{G}, d_{\max}) \to (\mathcal{G}, d_{\mathcal{G}})
$$
is continuous but not uniformly continuous.
\end{theorem}

An important question is how to understand the completion of 
$$
(\mathcal{HGM}, d_{\mathcal{HGM}}).
$$
To understand this question, we proposed a program for geometric Gibbs measures in~\cite{JGibbs0,JGibbsIntro,JGibbs}. However, we left the uniqueness conjecture in~\cite{JGibbs}. In this paper, we solve this conjecture and, consequently, we complete our project on the Teichmüller complex manifold structure on the space of all geometric Gibbs measures.    
We first briefly review geometric Gibbs measures, with more detail in~\cite{JGibbs}.

Recall that $\Si^{+}$ is a symbolic representation of $\T$.
The symbolic dynamical system $\sigma^{+}: \Si^{+}\to \Si^{+}$ is the symbolic representation of all circle endomorphisms of degree $d$ as dynamical systems.  
Let ${\mathcal M}^{+} (\sigma^{+})$ denote the space of all $\sigma^{+}$-invariant probability measures with $\mu_{0}^{+}$ as the basepoint.
The space $\Sigma$ and the shift map $\sigma$ serve as the dual symbolic representation of all circle endomorphisms of degree $d$. Let ${\mathcal M} (\sigma)$ denote 
the space of all $\sigma$-invariant probability measures with basepoint $\mu_{0}$. 

A measure $\mu \in {\mathcal M}(\sigma)$ (respectively, $\mu^{+}\in {\mathcal M}(\sigma^{+})$) is called non-atomic if $\mu (\{ v\} )=0$ for every point $v\in \Si$ (respectively, $\mu^{+}(\{w\})=0$ for every point $w\in \Si^{+})$. We say a measure $\mu \in {\mathcal M}(\sigma)$ 
(respectively, $\mu^{+}\in {\mathcal M}(\sigma^{+})$) has full support if $\mu ([v]_{n})\not=0$ for 
every $n$-cylinder $[v]_{n}$ in $\Si$ (respectively, $\mu^{+}([w]_{n})\not=0$ for every $n$-cylinder $[w]_{n}$ in $\Si^{+})$.  
 
For any $v=j_{0}j_{1}\cdots j_{n-1}\cdots\in \Si^{+}$, let $x=x(v)\in [0,1]$ in (\ref{rep}) and $z=e^{2\pi i x}\in \T$.
Consider a cylinder $[v]_{n}^{+}$ in $\Si^{+}$. For $v_{n}=j_{0}j_{1}\cdots j_{n-1}$, define  
$w_{n}=i_{n-1}\cdots i_{1}i_{0}$ where $i_{0}=j_{n-1}, \cdots, i_{n-1}=j_{0}$. Then we have the corresponding $n$-cylinder $[w]_{n}$ in $\Si$ where $w= \cdots w_{n}\in \Sigma$. 

Suppose $\mu\in {\mathcal M}(\sigma)$ is non-atomic and has full support.
For each $n\geq 0$, consider $d^{n+1}$ intervals labeled $J_{\underbrace{00...0}_{n}}$, $\cdots$, $J_{v_{n}}$, $\cdots$, $J_{\underbrace{(d-1)(d-1)...(d-1)}_{n}}$, where each interval has length
\begin{equation}~\label{length}
|J_{v_{n}}|=|J_{w_{n}}|=\mu ([w]_{n}).
\end{equation}
These intervals are arranged counter-clockwise in numerical order of the angle of $z\in J_{v_{n}}$ on the unit circle $\T$ beginning at $1$ (see~\cite[Figure 3]{JGibbs}).  
Since $\mu$ is $\sigma$-invariant, we have 
$$
|J_{v_{n}}| = \mu ([w]_{n}) =\sum_{j=0}^{d-1} \mu ([wj]_{n+1})=\cup_{j=0}^{d-1} |J_{v_{n}j}|,
$$
which implies that 
\begin{equation}~\label{lengthadd}
J_{v_{n}} = \cup_{j=0}^{d-1} J_{v_{n}j}
\end{equation}
and that the intervals nest as follows:
$$
\cdots \subset J_{v_{n}} \subset J_{v_{n-1}} \subset \cdots \subset J_{v_{1}}. 
$$
Since $\mu$ is non-atomic and has full support, the intersection $\cap_{k=1}^{\infty} J_{v_{n}}$ contains only a single point, which we denote as $h_{\mu}(z)$. This process defines a circle homeomorphism 
$$
h_{\mu}(z): \T\to \T.
$$ 
Then the map 
\begin{equation}~\label{fmu}
f_{\mu}=h_{\mu} \circ q_{d}\circ h_{\mu}^{-1}
\end{equation}
becomes a circle endomorphism of degree $d$ that preserves the Lebesgue measure $|\cdot|$ on $\T$, meaning
$$
|f^{-1}_{\mu}(J)| =|J|
$$
for every interval $J\subset \T$.

The condition that $f_{\mu}$ is uniformly symmetric in Definition~\ref{dggm} follows from the asymptotically symmetric balanced geometry of $\mu$ (see~\cite {JFurs}). First we define the adding machine $\hbox{add}: \Sigma_{n}\to \Sigma_{n}$ on the set of all $n$-cylinders
$
\Sigma_{n}=\{ [w_{n}]=[i_{n-1} \cdots i_{1}i_{0}]\},
$
by
$$
\hbox{add}([w_{n}])=
\begin{cases}
 [i_{n-1}\cdots i_{1} (i_{0}+1)], & \hbox{if $0\leq i_{0}<d-1$};\\
 [i_{n-1}\cdots (i_{k+1}+1)0\cdots 0], & \hbox{if $i_{k}=d-1$ for $0\leq k< n-1$};\\
 [0\cdots 0], &\hbox{if $i_{k}=d-1$ for all $0\leq k\leq n-1$} 
\end{cases}
$$
We say that $\mu$ has asymptotically symmetric balanced geometry if there is a bounded function $\epsilon (n) >0$ and $\epsilon (n)\to 0$ as $n\to \infty$ such that 
$$
\frac{1}{1+\epsilon (n)}\leq \frac{\mu (\hbox{add} ([w_{n}]))}{\mu ([w_{n}]} \leq 1+\epsilon (n)
$$
for all $w_{n}\in \Sigma_{n}$ and all $n\geq 1$. 

\medskip
\begin{proposition}~\label{asbg}
The circle endomorphism $f_{\mu}$ in (\ref{fmu}) is uniformly symmetric if and only if $\mu$ has asymptotically symmetric balanced geometry. 
\end{proposition}

\begin{proof}
Recall that for $w_{n}=i_{n-1}\cdots i_{1}i_{0}$, $v_{n}=j_{0}\cdots j_{n-2}j_{n-1}$ where $j_{n-1}=i_{0}$, $\cdots$, $j_{n-2}=i_{1}$, $j_{n-1}=i_{0}$. Suppose $\hbox{add}(w_{n})=i_{n-1}'\cdots i_{1}'i_{0}'$ and $v_{n}'=j_{0}'\cdots j_{n-2}'j_{n-1}'$ where $j_{0}'=i_{n-1}'$, $\cdots$, $j_{n-2}'=i_{1}'$, $j_{n-1}'=i_{0}'$. Then $J_{v_{n}}=J_{w_{n}}$ and $J_{v_{n}'}=J_{\hbox{add}(w_{n})}$ are two adjacent intervals in $\T$. Then that $\mu$ has asymptotically symmetric balanced geometry is equivalent to say that the sequence of Markov partitions $\{ \varpi_{n}=\{J_{w_{n}}\}\}$ has asymptotically symmetric bounded-nearby geometry (refer to~\cite{JBook,JGibbs}). From the proof that quasisymmetry on $f$ is equivalent to the sequence of Markov partitions having bounded-nearby geometry (refer to~\cite{JBook,JGibbs}), we see that all inverse branches $f^{-n}_{\mu}$ are uniformly symmetric is equivalent to the sequence of Markov partitions having asymptotically symmetric bounded-nearby geometry, which implies that $f_{\mu}$ is uniformly symmetric. 
\end{proof}

\medskip
\begin{definition}~\label{dggm}
Suppose $g\in {\mathcal G}$. A non-atomic $\sigma$-invariant probability measure $\mu\in {\mathcal M}(\sigma)$ is called a geometric Gibbs measure associated with $g$ if 
\begin{equation}~\label{eqggm}
\lim_{n\to \infty} \frac{\mu ([w]_{n})}{ \mu ([\sigma(w)]_{n-1})} = g(w) 
\end{equation}
uniformly on $w\in \Sigma$ and if $\mu$ has asymptotically symmetric balanced geometry.
Furthermore, if $f_{\mu}$ in (\ref{fmu}) is $C^{1+\alpha}$ for some $0<\alpha\leq 1$ and (\ref{eqggm}) holds, then we call $\mu$ a smooth geometric Gibbs measure associated with $g$.  
\end{definition}

\medskip
\begin{remark}
Note that when $\mu$ is a geometric Gibbs measure, then $h_{\mu}$ is quasisymmetric (see~\cite[Corollary 5.6]{JGibbs}). 
Moreover, if $f_{\mu}$ in (\ref{fmu}) is of class $C^{1+\alpha}$, then $f_{\mu}$ is expanding. Indeed, this follows from~\cite{Jkc}, since the conjugacy $h_{\mu}$ is H\"older continuous.    
\end{remark}

\medskip
\begin{remark}~\label{explain}
Technically, a geometric Gibbs measure $\mu$ associated with a continuous function $g\in{\mathcal G}$ is, in general, not a Gibbs state in the classical thermodynamic sense. That is, the measure of an $n$-cylinder need not be comparable to the exponential of the Birkhoff sum of the potential $\log g$ along the first $n$ iterates of an arbitrary point in the cylinder. However, when the limit in~(\ref{eqggm}) converges exponentially, the corresponding geometric Gibbs measure is indeed a Gibbs state in the classical sense.

Conceptually, a geometric Gibbs measure $\mu$ should be viewed instead as a finite geometric deformation of a reference measure on the underlying fractal set, whose infinitesimal scaling is prescribed by the function $g$.

Although a general geometric Gibbs measure $\mu$ need not satisfy the classical Gibbs property, it shares with the smooth case a fundamental geometric feature: when $\mu$ is regarded as defining a metric on the underlying fractal set, the associated dynamical system has derivative prescribed by the function $g$. This geometric characterization, rather than the classical Gibbs property, motivates the term \emph{geometric Gibbs measure} introduced in~\cite{JGibbs}.
\end{remark}

\section{\bf Existence and Uniqueness}~\label{uniq} 

\cite[Theorem 10.2]{JGibbs} proves the existence of a geometric Gibbs measure $\mu$ associated with each $g\in\mathcal{G}$. Moreover, in the same theorem, we also prove that for the basepoint $g_{0} =1/d$ in $\mathcal{G}$, there is only one geometric Gibbs measure associated with it (see also~\cite{JSymm}). Furthermore, if $g\in \mathcal{HG}$, it has a unique smooth geometric Gibbs measure, which is also the Gibbs measure in the classical sense; however, it may still have a geometric Gibbs measure that is not smooth.  One of the main results in this paper is the uniqueness of a geometric Gibbs measure for any given $g\in\mathcal{G}$.

\medskip
\begin{starthma}~\label{main}
For every $g\in\mathcal{G}$, there is a unique geometric Gibbs measure $\mu$ associated with $g$. 
\end{starthma}

\begin{proof}
\cite[Theorem 10.2]{JGibbs} establishes the existence of a geometric Gibbs measure associated with $g$. Assume $\mu_{1}$ and $\mu_{2}$ are both geometric Gibbs measures associated with $g$. 
We have that $f_{\mu_{1}}$ and $f_{\mu_{2}}$ belong to $\mathcal{US}$ and both preserve the Lebesgue measure, meaning that for any interval $J\subset \T$,
$$
|f_{\mu_{1}}^{-1}(J)| =|J|\quad \hbox{and}\quad |f_{\mu_{2}}^{-1} (J)|=|J|.
$$
Since both of $f_{\mu_{1}}$ and $f_{\mu_{2}}$ are circle endomorphisms of degree $d$ in $\mathcal{US}$, we have an orientation-preserving circle homeomorphism $h$ of $\T$ such that
$$
h\circ f_{\mu_{1}} =f_{\mu_{2}}\circ h.
$$
Furthermore, $h$ is quasisymmetric due to~\cite[Corollary 5.6]{JGibbs}.

The first part of the proof is to show that $h$ is symmetric in the sense of Definition~\ref{sch}. We state it as a lemma.

\medskip
\begin{lemma}~\label{symmclaim}
The homeomorphism $h$ is symmetric.
\end{lemma}

\begin{proof}

Let $F: \R\to \R$ be the lift of $f_{\mu_{1}}$ (or $f_{\mu_{2}}$). 
Let 
$$
\lambda = \frac{1}{g (\cdots 000)} =\lim_{n\to \infty}
\frac{F^{-n}(1)}{F^{-(n+1)}(1)}.
$$
The limit exists due to the quasisymmetric distortion lemma (see~\cite[Lemma 6.2 ]{JGibbs}) and the 
uniform symmetry of the sequence $\{F^{-n}\}_{n=1}^{\infty}$. 
Consider the linear map
$$
P(x) =\lambda x: \R\to \R.
$$
For each $n>0$, consider the homeomorphism.
$$
\gamma_{n}(x) =\frac{F^{-n}(x)}{F^{-n}(1)} : {\mathbb R}\to
{\mathbb R}.
$$
For any $\epsilon>0$, there is an $n_{0}>0$ such that $F^{-m}|[0, F^{-n}(1)]$
is $(1+\epsilon)$-quasisymmetric for any $m>n\geq n_{0}$. Then
$$
H(y) = \frac{F^{-m+n} (F^{-n}(1)y)}{F^{-m}(1)} :[0,1]\to[0,1]
$$
is a $(1+\epsilon)$-quasisymmetric homeomorphism with $H(0)=0$ and
$H(1)=1$. By the quasisymmetric distortion lemma (see~\cite[Lemma 6.2]{JGibbs}),  
$$
|H(y)-y| \leq \epsilon, \quad \forall y\in [0,1].
$$
Thus, for $y=F^{-n}(x)/F^{-n}(1)$, $x\in [0,1]$, we have
$$
|\gamma_{n+m}(x) - \gamma_{n}(x)|< \epsilon, \;\;\forall
m>n\geq n_{0}.
$$
This shows that the sequence $\{ \gamma_{n}(x)\}_{n=0}^{\infty}$ is
uniformly Cauchy on $[0,1]$ and, thus, converges
uniformly to a function $\gamma (x)$ on $[0,1]$. Similarly the
sequence of inverses $\{ \gamma_{n}^{-1}(y) = F^{n}(F^{-n}(1)
y)\}_{n=0}^{\infty}$ is also uniformly Cauchy on $[0,1]$ and converges uniformly to a function on $[0,1]$ which serves as the
inverse of $\gamma (x)$. So $\gamma (x)$ is a homeomorphism.
A direct calculation shows that $\gamma (x)$ is symmetric on $[0,1]$ since the sequence $\{ \gamma_{n}(x)\}_{n=1}^{\infty}$ is uniformly symmetric.

For any fixed $k>0$, $F^{-k}$ maps $[0, d^{k}]$ onto
$[0,1]$. Utilizing the relationship
$$
\gamma_{n+k} (x) = \frac{F^{-n}(1)}{F^{-n-k}(1)} \gamma_{n}
(F^{-k}(x)), \; x\in [0, d^{k}],
$$
we conclude that the sequence $\{ \gamma_{n}(x)\}_{n=0}^{\infty}$ is uniformly
Cauchy on $[0,d^{k}]$ and converges to a symmetric homeomorphism $\gamma (x)$ on $[0, d^{k}]$, given that $\{
F^{-k}\}_{k=0}^{\infty}$ is uniformly symmetric. Thus, $\{ \gamma_{n}(x)\}$ converges on any compact set in ${\mathbb R}^{+}=[0, \infty)$ and similarly in ${\mathbb R}^{-}=(-\infty, 0]$. It can be verified that 
$$
\gamma \circ F \circ \gamma^{-1 } (x) = P(x), \;\;
\forall x\in \R
$$
which says that $\gamma$ conjugates $F$ to the linear dynamical system $P(x)=\lambda x$.

Let 
\begin{equation}~\label{linearmodel}
L (x) =\gamma (\gamma^{-1} (x)+1).
\end{equation}
We have $L(0)=1$. The map $L$ is a symmetric homomorphism and 
\begin{equation}~\label{multiequation}
P ( L (x))  = L^{d} ( P (x)).
 \end{equation}

For any $w=\cdots i_{n-1}\cdots i_{1}i_{0}\in \Sigma$, let
$w_{n}=i_{n-1}\cdots i_{1}i_{0}$ and $\sigma(w_{n})=i_{n-1}\cdots i_{1}$. 
Since $\gamma(x)$ is symmetric on
$[0,1]$ with $\gamma(0)=0$ and $\gamma(1)=1$, from the definition of $\mu_{1}^{+}$ and $f_{\mu_{1}}$ (or $\mu_{2}^{+}$ and $f_{\mu_{2}}$), we have that
$$
g (w) = \lim_{n\to \infty}
\frac{|J_{w_{n}}|}{|J_{\sigma(w_{n}}|} =\lim_{n\to\infty}
\frac{|\gamma(I_{w_{n}})|}{|\gamma(I_{\sigma(w_{n}})|},
$$
where $I_{w_{n}}$ and $I_{\sigma(w_{n})}$ are the lifts of $J_{w_{n}}$ and $J_{\sigma(w_{n})}$ to $[0,1]$. 
Consider the non-negative integers $k$ and $k'$:
$$
k=i_{0}+i_{1}d+\cdots+i_{n-1}d^{n-1} \quad \hbox{and}\quad
k'=i_{1}+i_{2}d+\cdots +i_{n-1}d^{n-2}.
$$
We have that $k=dk'+i_{0}$ and that 
$$
I_{w_{n}} = F^{-n} ([k, k+1])\quad \hbox{and}\quad I_{\sigma(w_{n})} =
F^{-(n-1)} ([k', k'+1]).
$$
Using the property 
$$
\gamma (F^{-n}(x)) =\lambda^{-n}\gamma (x)
$$ 
and the fact that 
$$
\lambda \gamma (k') = \gamma (F(k')) =\gamma (dk'),
$$
we get 
$$
\frac{|\gamma(I_{w_{n}})|}{|\gamma (I_{\sigma(w_{n})})|}
=\frac{\lambda |\gamma ([k,k+1])|}{|\gamma ([k',k'+1])|} =
\frac{|\gamma ([k,k+1])|}{|\gamma ([k-i_{0},k+d -i_{0}])|}.
$$
This implies that
\begin{equation}~\label{lmdd}
\frac{|\gamma ([k,k+1])|}{|\gamma ([k-i_{0},k+d -i_{0}])|} =
g(\cdots 000i_{n-1}\cdots i_{1}i_{0}).
\end{equation}
Since $L^{k}(0)=\gamma (k)$, $g$ is completely determined by $\{ L^{k}(0)\}_{k=0}^{\infty}$ and vice versa since 
the set of points $\{ \cdots 000 i_{n-1}\cdots i_{1}i_{0} \}$ is dense in $\Si$. 
 
From (\ref{lmdd}), we have
$$
g (\cdots 000i_{n-1}\cdots i_{1}i_{0})= \frac{L^{k+1}(0)-L^{k}(0)}{L^{k+d-i_{0}}(0) -L^{k-i_{0}}(0)}
$$
This leads us to the conclusion that $L(0)=1$ and $L^{d} (0)=g_{f}(\cdots 000)^{-1}=\lambda$. Furthermore, we use  
$$
g (\cdots 00i) =\frac{L^{i+1} (0)-L^{i} (0)}{\lambda}, \quad 1\leq i\leq d-1,
$$
to show that $L^{i+1}_{f}(0)$ are uniquely determined by $g_{f}(\cdots 00i)$, $1\leq i\leq d-1$.
A similar reasoning extends to show that all $L^{k}(0)$, $1\leq k<\infty$ are uniquely determined 
by $\{g (\cdots 000i_{n-1}\cdots i_{1}i_{0})\}$, inductively.  

From Equation (\ref{multiequation}), we have that
$\lambda^{n} = L^{d^{n}} (0)$. We derive a set ${\mathcal A}=\{\lambda^{-n} L^{k}(0)\}_{k=0}^{d^{n}}$ giving a partition of $[0,1]$ and another set ${\mathcal B}=\{\lambda^{-(n-1)} L^{k}(0)\}_{k=0}^{d^{n}}$ giving a partition of $[0, \lambda]$. 
Consequently, $L$ maps the partition points of ${\mathcal A}$ to those of ${\mathcal B}$. 
Inductively, as established, we see $L$ is uniquely determined by $g$.  

Since $\mu_{1}$ and $\mu_{2}$ represent geometric Gibbs measures associated with the same potential $g$, they  share the same $L$. Let
$\gamma_{\mu_{1}}$ and $\gamma_{\mu_{2}}$ 
be the symmetric homeomorphisms constructed from $f_{\mu_{1}}$ and $f_{\mu_{2}}$, respectively. We then have the relationship: 
$$
L (x) =\gamma_{\mu_{1}} (\gamma_{\mu_{1}}^{-1} (x)+1)=\gamma_{\mu_{2}} (\gamma_{\mu_{2}}^{-1} (x)+1).
$$
This indicates that the lift $H$ of $h$ is
$$
H =\gamma_{\mu_{2}}^{-1}\circ \gamma_{\mu_{1}},
$$
which is symmetric since both $\gamma_{\mu_{1}}$ and $\gamma_{\mu_{2}}$ are symmetric (see Remark~\ref{cmp}), thereby confirming that 
$h$ is symmetric. We proved the lemma. 
\end{proof}

The second part of the proof is to prove that $h$ is the identity. We state it as a lemma.

\medskip
\begin{lemma}~\label{idclaim}
The homeomorphism $h$ is the identity.
\end{lemma}

\begin{proof}
To demonstrate that $h$ is the identity using the symmetry property for $h$ and the uniform symmetry property for $f_{\mu_{1}}$ and $f_{\mu_{2}}$, we start with the definition of $\Phi$:
$$
1\leq \Phi=\sup_{I\subseteq \T} \frac{|h(I)|}{|I|}\leq \infty.
$$
The set $X$ is defined as:
$$
 X=\{ z=e^{2\pi i x}\in \T\;|\; \exists I_{k}^{z}=\arc{e^{2\pi i a_{k}}e^{2\pi i b_{k}}}, \lim_{k\to \infty} a_{k}=\lim_{k\to \infty} b_{k} =x, \lim_{k\to \infty} \frac{|h(I_{k}^z)|}{|I_{k}^z|} =\Phi\}.
 $$

Given the sequence of the partitions $\{\eta_{n}\}_{n=1}^{\infty}$ for $f_{\mu_{1}}$, $h(\eta_{n})$ is the corresponding sequence of the partitions for $f_{\mu_{2}}$ (see Section 1). The uniform symmetry of $f_{\mu_{1}}$ and $f_{\mu_{2}}$ implies the bounded geometry property for both measures, characterized by a constant $C>0$ such that for $J\subset I$, $J\in \eta_{n, f_{\mu_{1}}}$ and $I\in \eta_{n-1, f_{\mu_{1}}}$ (or $J\in \eta_{n, f_{\mu_{2}}}$ and $I\in \eta_{n-1, f_{\mu_{2}}}$), $n\geq 1$, 
$$ 
 \frac{|I|}{|J|} \leq C.
 $$  
 This bounded geometry ensures that for all $w=\cdots w_{n}\in \Si$ and all $n\geq 1$,
 $$
 \frac{|I_{v_{n}}|}{|I_{\s(v_{n})}|} \geq \frac{1}{C}\quad \hbox{and} \quad  \frac{|h(I_{v_{n}})|}{|h(I_{\s(v_{n})})|}\geq \frac{1}{C}.  
 $$
 In this context, recall that when $w_{n} =i_{n-1}i_{n-2}\cdots i_{1}i_{0}$,  $v_{n}$ is defined as $v_{n} =j_{0}j_{1}\cdots j_{n-2}j_{n-1}$ for $j_{0}=i_{n-1}$, $j_{1}=i_{n-2}$, $\cdots$, $j_{n-2}=i_{1}$, $j_{n-1}=i_{0}$. Then $\s (v_{n})= j_{0}j_{1}\cdots j_{n-2}$ and we have $I_{v_{n}}\subset I_{\s (v_{n})}$.

Suppose $\{I_{k}=\arc{e^{2\pi i a_{k}}e^{2\pi i b_{k}}}\}_{k=1}^{\infty} $ is a sequence of arcs such that 
$$
\lim_{k\to \infty} \frac{|h(I_{k})|}{|I_{k}|}=\Phi.
$$
We assume, by taking a subsequence if necessary, that $\{ a_{k}\}_{k=1}^{\infty}$ and  $\{b_{k}\}_{k=1}^{\infty}$ converge to $a$ and $b$, respectively. 
If $a=b=x$, then $z=e^{2\pi i x}\in X$. (Note that if $\Phi=\infty$ then $a=b$.) 
If $a<b$, $I=\arc{e^{2\pi i a}e^{2\pi i b}}$ is a non-trivial arc satisfying 
$$
\frac{|h(I)|}{|I|}=\Phi,
$$
and consequently, all interior points of $I$ are in $X$, proving that $X$ is non-empty. 
 
To establish that $X$ is dense in $\T$, we claim that for any $I_{v_{n}}\in \eta_{n}$, $X\cap I_{v_{n}} \not= \emptyset$.  
Before proving this claim, we first prove that   
\begin{equation}~\label{zero1}
\lim_{k\to \infty} \sum_{v_{n}^{1} \not=v_{n}} \cdots \sum_{v_{n}^{k}\not=v_{n}} |I_{v_{n}^{1}\cdots v_{n}^{k}}|=0,
\end{equation}
where the $v_n^i$ are all words of length $n$.
First, we have that 
 $$
 |I_{v_{n}}| \geq A=\frac{1}{C^{n}}.
 $$
 Since 
 $$
 \bigcup_{v_{n}^{1}\not=v_{n}} I_{v_{n}^{1}}=\overline{\T\setminus I_{v_{n}}},
 $$  
 we get
 $$
 \sum_{v_{n}^{1}\not=v_{n}} |I_{v_{n}^{1}}| = 1-|I_{v_{n}}|  \leq 1-A.
 $$
For any $v_{n}^{1}$, we have that $I_{v_{n}^{1}v_{n}} \subset I_{v_{n}^{1}}$ and  
$$
|I_{v_{n}^{1}v_{n}}| \geq A |I_{v_{n}^{1}}|.
$$ 
Since 
$$
\bigcup_{v_{n}^{2}\not=v_{n}} I_{v_{n}^{1}v_{n}^{2}}=\overline{I_{v_{n}^{1}}\setminus I_{v_{n}^{1}v_{n}}},
$$
we have
$$
\sum_{v_{n}^{2}\not=v_{n}} |I_{v_{n}^{1}v_{n}^{2}}|= |I_{v_{n}^{1}}|-|I_{v_{n}^{1}v_{n}}| \leq |I_{v_n^1}| - A|I_{v_n^1}| = (1-A) |I_{v_{n}^{1}}|, 
$$
which implies that 
$$
\sum_{v_{n}^{1}\not=v_{n}} \sum_{v_{n}^{2}\not= v_{n}} |I_{v_{n}^{1}v_{n}^{2}}|\leq (1-A)\sum_{v_{n}^{1}\not=v_{n}} |I_{v_{n}^{1}}|\leq (1-A)^{2}.
$$ 
Inductively, suppose we have
$$
\sum_{v_{n}^{1}\not=v_{n}} \cdots \sum_{v_{n}^{k-1}\not=v_{n}} |I_{v_{n}^{1}\ldots v_{n}^{k-1}}| 
\leq (1-A)^{k-1}
$$
for $k\geq 3$.
Then 
$$
\sum_{v_{n}^{1}\not=v_{n}} \cdots \sum_{v_{n}^{k}\not=v_{n}} |I_{v_{n}^{1}\ldots v_{n}^{k}}|
 =\sum_{v_{n}^{1}\not=v_{n}} \cdots \sum_{v_{n}^{k-1}\not=v_{n}} \Big( |I_{v_{n}^{1}\cdots v_{n}^{k-1}}| - |I_{v_{n}^{1}\cdots v_{n}^{k-1}v_{n}}|\Big)
$$
$$
= \sum_{v_{n}^{1}\not=v_{n}} \cdots \sum_{v_{n}^{k-1}\not=v_{n}} |I_{v_{n}^{1}\cdots v_{n}^{k-1}}| \Big( 1-\frac{  |I_{v_{n}^{1}\cdots v_{n}^{k-1}v_{n}}|}{  |I_{v_{n}^{1}\cdots v_{n}^{k-1}}|}\Big).
$$
Notice that $I_{v_{n}^{1}\cdots v_{n}^{k-1}v_{n}} \subset I_{v_{n}^{1}\cdots v_{n}^{k-1}}$. Since 
$$
\frac{ |I_{v_{n}^{1}\cdots v_{n}^{k-1}v_{n}}|}{ |I_{v_{n}^{1}\cdots v_{n}^{k-1}}|}\geq A,
$$
$$
1-\frac{ |I_{v_{n}^{1}\cdots v_{n}^{k-1}v_{n}}|}{ |I_{v_{n}^{1}\cdots v_{n}^{k-1}}|}\leq 1-A.
$$
Thus
$$
\sum_{v_{n}^{1}\not=v_{n}} \cdots \sum_{v_{n}^{k}\not=v_{n}} |I_{v_{n}^{1}\ldots v_{n}^{k}}| \leq (1-A) \sum_{v_{n}^{1}\not=v_{n}} \cdots \sum_{v_{n}^{k-1}\not=v_{n}} |I_{v_{n}^{1}\cdots v_{n}^{k-1}}|\leq (1-A)^{k}.
$$
Letting $k\to\infty$, we have that 
 $$
\lim_{k\to \infty} \sum_{v_{n}^{1} \not=v_{n}} \cdots \sum_{v_{n}^{k}\not=v_{n}} |I_{v_{n}^{1}\cdots v_{n}^{k}}|=0.
$$ 

Now we are ready to prove the claim by contradiction. Suppose $I_{v_{n}}\cap X=\emptyset$. 
Then we can find a number $D<\Phi$ such that 
$$
\frac{|h(I)|}{|I|} \leq D, \quad \forall \; I\subset I_{v_{n}}.
$$ 
Since $X\not=\emptyset$, we have an interval $I^{D}\subseteq \T$
such that 
$$
\frac{|h(I^{D})|}{|I^{D}|} > D.
$$
We pull back $I^{D}$ by $f^{n}_{\mu_{1}}$ to get
$f^{-n}_{\mu_{1}} (I^{D})=\cup_{v_{n}^{1}} I_{v_{n}^{1}}^{D}$, where $I^{D}_{v_{n}^{1}}\subset I_{v_{n}^{1}}\in \eta_{n}$ and $f^{n}_{\mu_{1}} (I^{D}_{v_{n}^{1}}) =I^{D}$.
This implies that for all $k\geq2$,
$$
|I^{D}| = |I_{v_n}^D| + \sum_{v_n^1 \neq v_n}|I_{v_n^1}^D| = |I_{v_n}^D| + \sum_{v_n^1 \neq v_n}|f^{-n}(I_{v_n^1}^D)| 
$$
$$
= |I_{v_n}^D| + \sum_{v_n^1 \neq v_n}\left( |I_{v_n v_n^1}^D| + \sum_{v_n^2 \neq v_n}|I_{v_n^2v_n^1}^D| \right)
$$
$$
 = |I_{v_n}^D| + \sum_{v_n^1 \neq v_n}|I_{v_n v_n^1}^D| + \sum_{v_n^2\neq v_n}\sum_{v_n^1\neq v_n}|I_{v_n^2 v_n^1}^D| = \ldots
 $$
 $$
 = |I_{v_{n}}^{D}| +  \sum_{l=1}^{k-1}  \sum_{v_{n}^{l}\not=v_{n}} \cdots \sum_{v_{n}^{1}\not= v_{n}}  |I_{v_{n}v_{n}^{l}\ldots v_{n}^{1}}^{D}| + \sum_{v_{n}^{k}\not=v_{n}} \cdots \sum_{v_{n}^{1}\not= v_{n}} |I^{D}_{v_{n}^{k}\ldots v_{n}^{1}}|
$$
and, similarly, 
$$
|h(I^{D})|
=|h(I_{v_{n}}^{D})| + \sum_{l=1}^{k-1}  \sum_{v_{n}^{l}\not=v_{n}} \cdots \sum_{v_{n}^{1}\not= v_{n}}  |h(I_{v_{n}v_{n}^{l}\ldots v_{n}^{1}}^{D})| 
$$
$$
+  \sum_{v_{n}^{k}\not=v_{n}} \cdots \sum_{v_{n}^{1}\not= v_{n}} |h(I^{D}_{v_{n}^{k}\ldots v_{n}^{1}})|.
$$
Because $I_{v_{n}}^{D}$ and $I_{v_{n}v_{n}^{l}\ldots v_{n}^{1}}^{D}$ are sub-intervals of $I_{v_{n}}$,  
$$
\frac{|h(I^{D}_{v_{n}})| }{|I^{D}_{v_{n}}|},  \;\;
\frac{|h(I^{D}_{v_{n}v_{n}^{l}\cdots v_{n}^{1}})| }{|I^{D}_{v_{n}v_{n}^{l}\cdots v_{n}^{1}}|} \leq D,\quad \forall \; l\geq 1,
$$
which implies that 
$$
\frac{|h(I_{v_{n}}^{D})| + \sum_{l=1}^{k-1}  \sum_{v_{n}^{l}\not=v_{n}} \cdots \sum_{v_{n}^{1}\not= v_{n}}  |h(I_{v_{n}v_{n}^{l}\ldots v_{n}^{1}}^{D})|}{|I_{v_{n}}^{D}| +  \sum_{l=1}^{k-1}  \sum_{v_{n}^{l}\not=v_{n}} \cdots \sum_{v_{n}^{1}\not= v_{n}}  |I_{v_{n}v_{n}^{l}\ldots v_{n}^{1}}^{D}|}\leq D, \quad \forall\; k\geq 2.
$$
Applying 
$$
\lim_{k\to \infty} \sum_{v_{n}^{k}\not=v_{n}} \cdots \sum_{v_{n}^{1}\not= v_{n}} |I^{D}_{v_{n}^{k}\cdots v_{n}^{1}}|=0
$$
and 
$$
\lim_{k\to \infty}  \sum_{v_{n}^{k}\not=v_{n}} \cdots \sum_{v_{n}^{1}\not= v_{n}} |h(I^{D}_{v_{n}^{k}\cdots v_{n}^{1}})|=0,
$$
 we get
$$
\frac{|h(I^{D})|}{|I^{D}|} 
\leq D.
$$
This is a contradiction. The contradiction completes the proof of the claim.  
 
Finally, we will prove $\Phi=1$. Again, we will prove this by using contradiction. 
Assume $\Phi>1$. Then we have two numbers $1< D_{1}< D_{2}<\Phi$.
Since $h$ is symmetric, there exists a positive bounded function $\epsilon(t)$
  such that $\epsilon(t)\to0$ as $t\to0^{+}$ and
$$
    \frac{1}{1+\epsilon(t)}\leq
    \frac{|h(I)|}{|h(I')|}\leq1+\epsilon(t)
    $$
holds for all closed arcs $I=\arc{e^{2\pi i (x-t)}e^{2\pi i x}}$ and $I'=\arc{e^{2\pi i x} e^{2\pi i (x+t)}}$.
Fix $t_0$ such that 
$$
  \epsilon(t)<\frac{D_2}{D_1}-1 \quad \forall t<t_0.
$$
Since $\overline{X}\not=\emptyset$, there exists an arc
$I=\arc{e^{2\pi i a}e^{2\pi i b}}\subset \T$ with $|I|<t_0$ such that 
$$
\frac{|h(I)|}{|I|} > D_{2}.
$$
Let $L=\arc{e^{2\pi i(2a-b)}e^{2\pi i a}}$ and $R=\arc{e^{2\pi i b}e^{2\pi i (2b-a)}}$. We have that 
$$
  \frac{|h(R)|}{|R|} = \frac{|h(R)|}{|h(I)|} \cdot \frac{|h(I)|}{|I|} \cdot \frac{|I|}{|R|} > \frac{1}{1+\epsilon(|I|)} \cdot D_2 \cdot 1 > D_1
$$
and
$$
  \frac{|h(L)|}{|L|} = \frac{|h(L)|}{|h(I)|} \cdot \frac{|h(I)|}{|I|} \cdot \frac{|I|}{|L|} > \frac{1}{1+\epsilon(|I|)} \cdot D_2 \cdot 1 > D_1.
$$

Now we want to show that for $I_{a} =\arc{e^{2\pi i a}1}$
\begin{equation*}
  \frac{|h(I_{a})|}{|I_{a}|} > D_1.
\end{equation*}
Consider any arc $J=\arc{e^{2\pi i b}e^{2\pi i c}}\supset R$ with $2b-a\leq c\leq 1$ satisfying  
$$
 \frac{|h(J)|}{|J|} >D_1.
$$
If $c=1$, we have 
$$
\frac{|h(I_{a})|}{|I_{a}|} = \frac{|h(I\cup J)|}{|I\cup J|} >D_{1}.
$$
Then we have nothing further to prove.
If $c<1$, we have a number $\delta >0$ such that $c+\delta< 1$ and such that for any $x\in [c, c+\delta]$ and $J_{x} =\arc{e^{2\pi i b}e^{2\pi x}}$, we have 
$$
 \frac{|h(J_{x})|}{|J_{x}|} >D_1.
$$
Since $\overline{X} =\T$, there is an arc $I_{1}=\arc{e^{2\pi i a_{1}}e^{2\pi b_{1}}} \subset \arc{e^{2\pi i c}e^{2\pi i (c+\delta)}}$ with $|I_1|<t_0$ such that 
\begin{equation*} 
 \frac{|h(I_{1})|}{|I_{1}|} >D_2.
\end{equation*}
Let $J_{1}=\arc{e^{2\pi i b}e^{2\pi i a_{1}}}$. Then we have three consecutive arcs $I$, $J_{1}$, and $I_{1}$ such that 
$$
 \frac{|h(I \cup J_{1}\cup I_{1})|}{|I\cup J_{1}\cup I_{1}|}> D_{1}.
$$
Consider $I_{1}$ as a new $I$ and repeat the above construction. We get three consecutive intervals 
$I_{1}$, $J_{2}$, and $I_{2}$ such that
 $$
 \frac{|h(I_{1} \cup J_{2}\cup I_{2})|}{|I_{1}\cup J_{2}\cup I_{2}|}> D_{1}.
$$
Inductively, for each integer $n\geq 2$, we have three consecutive intervals 
$I_{n-1}$, $J_{n}$, and $I_{n}$ such that 
$$
 \frac{|h(I_{n-1} \cup J_{n}\cup I_{n})|}{|I_{n-1}\cup J_{n}\cup I_{n}|}> D_{1},
$$
which implies that 
$$
 \frac{|h(I\cup (\cup_{i=1}^{n} (J_{i}\cup I_{i})))|}{|I\cup (\cup_{i=1}^{n} (J_{i}\cup I_{i}))|}> D_{1}.
$$
If $I_{n}$ has an endpoint $1$ for some $n>1$, then  
we have 
$$
\frac{|h(I_{a})|}{|I_{a}|} >D_{1}.
$$
We have nothing further to prove.
If all endpoints of $I_{n}$ are not $1$ for all $n\geq 1$, suppose $I_{n}=\arc{z_{n}'z_{n}}$.
since $\{z_{n}<1\}_{n=1}^{\infty}$ is a strictly increasing sequence in $\T$, we have
$$
z_{\infty} =\lim_{n\to \infty} z_{n}\leq 1
$$
and 
$$
 \frac{|h(I\cup (\cup_{n=1}^{\infty} (J_{n}\cup I_{n})))|}{|I\cup (\cup_{n=1}^{\infty} (J_{n}\cup I_{n}))|}> D_{1}.
$$ 
Here we use the fact that if $\alpha_{1}, \alpha_{2}, \beta_{1}, \beta_{2}$ are positive numbers such that $\alpha_{1}/\alpha_{2}> \gamma$ and $\beta_{1}/\beta_{2}\geq \gamma$, then $(\alpha_{1}+\alpha_{2})/(\beta_{1}+\beta_{2}) >\gamma$. 
Since $z_{\infty}$ depends on the initially chosen arc $J$, we write it as $z_{\infty}(J)$.
Consider the set 
$$
{\mathcal B} =\{ z_{\infty} (J)\;|\; J \hbox{ satisfies the last inequality}\}
$$
Let $w =\sup {\mathcal B}$. We claim $w =1$. Otherwise, we take $J=\arc{e^{2\pi i b}e^{2\pi i w}}$. It satisfies the last inequality. Then $z_{\infty} (J) >w$. 
This is a contradiction. So 
$$
\frac{|h(I_{a})|}{|I_{a}|}  > D_{1}.
$$

Similarly, by using $L$ instead of $R$ and applying the procedure above, we get 
$$
\frac{|h(\T\setminus I_{a})|}{|\T\setminus I_{a}|}  > D_{1}.
$$
Finally, we get the following contradiction.
$$
1= \frac{|h(\T)|}{|\T|} =\frac{|h(\T\setminus I_{a})|+|h(I_{a})|}{|\T\setminus I_{a}|+|I_{a}|} >D_{1}>1.
$$
The contradiction implies that $\Phi=1$. 

Since $\Phi=1$, we have that for any non-trivial arc $J\subset \T$, $|h(J)|/|J| = 1$. Otherwise, if there is an interval $J\subset \T$ such that $|h(J)|/|J| <1$, let $L\cup R=\T \setminus J$. 
Then 
$$
1=\frac{|h(\T)|}{|\T|} = \frac{|h(L)|+|h(J)|+|h(R)|}{|L|+|J|+|R|} <1,
$$
since $|h(L)|/|L|\leq 1$ and $|h(R)|/|R|\leq 1$. This is a contradiction. Since $h(1)=1$, it follows that $h=id$. We proved the lemma.
\end{proof}

Lemma~\ref{symmclaim} and Lemma~\ref{idclaim} imply that $f_{\mu_{1}}=f_{\mu_{2}}$ and $h_{\mu_{1}} =h_{\mu_{2}}$. Thus $\mu_{1}=\mu_{2}$.  We completed the proof of the theorem.
\end{proof}

\begin{remark}
Lemma~\ref{idclaim} is a symmetric rigidity result. The first symmetric rigidity result was proved in~\cite{JSymm} when one of $f_{\mu_{1}}$ and $f_{\mu_{2}}$ is the polynomial $q_{d}(z)=z^{d}$, which implies that the geometric measure associated with the constant $g$-function $ g\equiv 1/d$ is unique. We proved a generalized version of the symmetric rigidity result in a joint work with my former Ph.D. students for circle endomorphisms with bounded geometry in~\cite{AHJW}. 
 \end{remark}

\section{\bf Complex Characterization of Uniform Symmetry}~\label{ccus}

The definition of a uniformly symmetric circle endomorphism (Definition~\ref{usce}) is in the sense of real analysis. In this section, we will give a complex characterization. For this purpose, we first define a quasisymmetric homeomorphism on an interval. 

\medskip
\begin{definition}~\label{intqs}
Suppose $Z: [0,1]\to [0,1]$ is a homeomorphism with $Z(0)=0$ and $Z(1)=1$. We say it is $M$-quasisymmetric for some $M\geq 1$ if 
$$
M^{-1} \leq \frac{|Z(x+t)-Z(x)|}{|Z(x)-Z(x-t)|} \leq M, \quad \forall x, x-t, x+t\in [0,1], t>0.
$$
\end{definition}

The proof of the following optimal quasisymmetric distortion result can be found in~\cite[Lemma 6.2]{JGibbs} or in~\cite{JSymm}.

\medskip
\begin{lemma}~\label{intqslem}
Suppose $Z: [0,1]\to [0,1]$ is a $M$-quasisymmetric homeomorphism with $Z(0)=0$ and $Z(1)=1$. Then 
$$
|Z(x)-x| \leq M-1, \quad \forall\; x\in [0,1].
$$
\end{lemma}

Suppose $h$ is a circle homeomorphism and $H$ is its lift. We denote
$$
\lambda_{H}(x, y, k) =\frac{|H(x+ky)-H(x)|}{|H(x)-H(x-y)|}, \quad x, y\in\R,\;0\leq k\leq 1,
$$
as the quasisymmetric distortion function of $H$. Then

\medskip
\begin{corollary}~\label{asymm1}
Suppose $h$ is an $M$-quasisymmetric circle homeomorphism. Then 
$$
| \lambda_{H} (x, y, k) -k|\leq M^{2}-1, \quad \forall x, y\in \R, \; 0\leq k\leq 1.
$$
\end{corollary}

\begin{proof}
Consider
$$
Z(k)= \frac{|H(x+ky)-H(x)|}{|H(x+y)-H(x)|}: [0,1]\to [0,1].
$$
It is a $M$-quasisymmetric homeomorphism of $[0,1]$ with $Z(0)=0$ and $Z(1)=1$. From Lemma~\ref{intqslem}, we have that 
$$
k -(M-1)\leq \frac{|H(x+ky)-H(x)|}{|H(x+y)-H(x)|} \leq k+(M-1).
$$
This implies that 
$$
\Big( k -(M-1)\Big) \lambda_{H} (x,y,1) \leq \frac{|H(x+ky)-H(x)|}{|H(x)-H(x-y)|} \leq \Big( k+(M-1)\Big) \lambda_{H}(x,y,1).
$$
Since 
$M^{-1} \leq \lambda_{H} (x, y,1) \leq M$,
$$
M^{-1}\Big( k -(M-1)\Big) \leq \frac{|H(x+ky)-H(x)|}{|H(x)-H(x-y)|} \leq M\Big( k+(M-1)\Big).
$$
Therefore, we have
$$
\Big|\lambda_{H} (x, y, k)-k\Big| \leq M^{2}-1.
$$
\end{proof}

\medskip
\begin{corollary} A circle homeomorphism $h$ is symmetric if and only if there is a bounded positive real function $\epsilon (t)$ with $\epsilon (t)\to 0^{+}$ as $t\to 0^{+}$ such that  
$$
| \lambda (x, y, k) -k|\leq \epsilon (|y|), \quad \forall x, y\in \R, \; 0\leq k\leq 1.
$$
\end{corollary}

Suppose $\Omega$ is a domain in the complex plane and $G: \Omega \to G(\Omega)$ is an orientation-preserving homeomorphism. We say it is quasiconformal if it has locally integrable distribution first partial derivatives $G_{z}$ and $G_{\overline{z}}$ satisfying the inequality 
$$
|G_{\overline{z}}(z)| \leq k |G_{z} (z)|
$$
for some $0\leq k <1$ and for almost all $z\in \Omega$. 
Let
$$
b_{G} (z) =\frac{G_{\overline{z}}(z)}{G_{z}(z)}
$$
be the Beltrami coefficient of $G$. The quasiconformal dilatation of $G$ is
$$
K_{G}(z) =\frac{1+|b_{G}(z)|}{1-|b_{G}(z)|}.
$$
Then we have
$$
k_{G}=\|b_{G}(z)\|_{\infty}\leq k,\quad K_{G} =\frac{1+k_{G}}{1-k_{G}} \leq K=\frac{1+k}{1-k}.
$$

Suppose $h$ is a quasisymmetric circle homeomorphism and $H$ is its lift. Let $\HH$ be a quasiconformal extension of $H$ to the upper-half plane $\UH$ with 
$$
\HH (z+1)=\HH (z)+1, \quad \forall z\in \UH.
$$
Then $\HH$ is the lift of a quasiconformal extension $\widehat{h}$ of $h$ to the unit disk $\D$. The lift $\HH$ of every quasiconformal extension $\widehat{h}$ of $h$ to $\D$ has this property. Among all quasiconformal extensions of $h$, we use 
$$
\widetilde{H} (z) = U (x,y)+iV (x,y), \quad z=x+iy\in \UH,
$$
where 
$$
U(x, y)= \frac{1}{2y} \int_{x-y}^{x+y} H(s) ds =\frac{1}{2} \int_{-1}^{1} H(x+ky) dk
$$
and
$$
V(x,y)=\frac{1}{y} \int_{x}^{x+y} H(s) ds -\frac{1}{y} \int_{x-y}^{x} H(s) ds,
$$
to denote the Beurling-Ahlfors extension. One can check:
$$
\widetilde{H}(z+1) =\widetilde{H} (z) +1, \quad \forall\; z\in \UH.
$$
Note that we use a slightly different normalization from that given in~\cite{AhlforsBook}. 
In our normalization, 
$$
\widetilde{Id_{\R}}=id_{\UH}, \quad 
\widetilde{A\circ H\circ B} =A\circ \widetilde{H}\circ B
$$
for any affine maps $A$ and $B$ on $\R$. 

For the quasisymmetric distortion function $\lambda_{H} (x, y, k)$, we have that
$$
\int_{0}^{1} \lambda_{H} (x,y,k) dk = \frac{1}{|H(x)-H(x-y)|}\Big| \frac{1}{y} \int_{x}^{x+y} H(s) ds- H(x)\Big|.
$$

To calculate the quasiconformal dilataion of $\widetilde{H}$, we have that
$$
K_{\widetilde{H}} (z) =\frac{\widetilde{H}_{z}+\widetilde{H}_{\overline{z}}}{\widetilde{H}_{z}-\widetilde{H}_{\overline{z}}}
=
\frac{(U+iV)_{z}+(U+iV)_{\overline{z}}}{ (U+iV)_{z}-(U+iV)_{\overline{z}}}
$$
$$
=\frac{(U+iV)_{x}- i(U+iV)_{y}+(U+iV)_{x}+
i(U+iV)_{y}}{(U+iV)_{x}- i(U+iV)_{y}-(U+iV)_{x}-
i(U+iV)_{y}}
$$
$$
=\frac{U_{x}+ iV_{x}}{V_{y}- iU_{y}}= \frac{1+ia}{b-ic}
$$
where $a=V_{x}/U_{x}$, $b=V_{y}/U_{x}$, and $c=U_{y}/U_{x}$. 

For $y>0$, let
$$
L= H(x)-H(x-y), \quad R=H(x+y)-H(x)
$$
and 
$$
L' =H(x)-\frac{1}{y} \int_{x-y}^{y} H(s) ds,\quad  R'= \frac{1}{y} \int_{x}^{x+y} H(s) ds -H(x).
$$
Then 
$$
U_{x}= \frac{1}{2y} (R+L), \quad V_{x} =\frac{1}{y} (R-L)
$$
and 
$$
V_{y} =\frac{1}{y} (R-L) -\frac{1}{y} (R'+L'), \quad U_{y} =\frac{1}{2y}(R-L)-\frac{1}{2y} (R'-L').
$$
Let 
$$
\rho =\lambda_{H} (x, y, 1) =\frac{R}{L},
$$
$$
\rho_{+} = \rho_{+}(x, y) =\frac{R'}{L},
$$
and
$$
\rho_{-} =\rho_{-} (x, y) = \frac{L'}{R}.
$$
Then we have
$$
a(1+\rho) =2\frac{R-L}{R+L} \frac{R+L}{L},
$$
$$
b(1+\rho) =2 \frac{R+L-R'-L'}{R+L} \frac{R+L}{L}
=2\Big(\frac{R}{L} +1-\frac{R'}{L} -\frac{R}{L}\frac{L'}{R}\Big),
$$
$$
c(1+\rho) =\frac{R-L-R'+L'}{R+L} \frac{R+L}{L} =\frac{R}{L}-1-\frac{R'}{L}+\frac{R}{L}\frac{L'}{R}.
$$
Finally, we get
$$
a=\frac{2(\rho-1)}{\rho+1},
$$
$$
b=\frac{2(\rho +1-\rho_{+} -\rho \rho_{-})}{ \rho+1},
$$
$$
c= \frac{\rho-1-\rho_{+}+\rho\rho_{-}}{\rho+1}.
$$
Since 
$$
K_{\widetilde{H}} (z) = \frac{1+ia}{b-ic},
$$
we have that
$$
K_{\widetilde{H}} (z) +1= \frac{1+ia+b-ic}{b-ic}.
$$

Suppose $h_{0}$ and $h_{1}$ are two quasisymmetric circle homeomorphisms and $H_{0}$ and $H_{1}$ are their lifts. Let 
$\widetilde{H}_{0}$ and $\widetilde{H}_{1}$ be their Beurling-Ahlfors extensions. Then the difference of their Beltrami coefficients, 
$$
b_{\widetilde{H}_{1}} (z)-b_{\widetilde{H}_{0}} (z)= \frac{K_{\widetilde{H}_{1}} (z)-1}{K_{\widetilde{H}_{1}} (z)+1}- \frac{K_{\widetilde{H}_{0}} (z)-1}{K_{\widetilde{H}_{0}} (z)+1}=2 \frac{K_{\widetilde{H}_{1}} (z)-K_{\widetilde{H}_{0}} (z)}{(K_{\widetilde{H}_{1}} (z)+1)(K_{\widetilde{H}_{0}} (z)+1)}
$$
$$
=2\frac{(1+ia_{1})(b_{0}-ic_{0}) -(1+ia_{0})(b_{1}-ic_{1})}{(1+ia_{1}+b_{1}-ic_{1})(1+ia_{0}+b_{0}-ic_{0})}
$$
$$
=2\frac{(a_{1}-a_{0})(ib_{1}+c_{1}) +(b_{0}-b_{1})(1+ia_{1})+(c_{1}-c_{0})(i-a_{1})}{(1+ia_{1}+b_{1}-ic_{1})(1+ia_{0}+b_{0}-ic_{0})}
$$
Since $R'_{i}\leq R_{i}$ and $L'_{i}\leq L_{i}$ for $i=0, 1$, we have that $\rho_{+,i} \leq \rho_{i}$ and $\rho_{i}\rho_{-,i}\leq 1$ for $i=0,1$. Thus $b_{0}, b_{1}>0$. This implies that
$$
|(1+ia_{1}+b_{1}-ic_{1})(1+ia_{0}+b_{0}-ic_{0})|\geq 1
$$
and
$$
|b_{\widetilde{H}_{1}} (z)-b_{\widetilde{H}_{0}} (z)|\leq 
2|(a_{1}-a_{0})(ib_{1}+c_{1}) +(b_{0}-b_{1})(1+ia_{1})+(c_{1}-c_{0})(i-a_{1})|
$$
$$
\leq C (|a_{1}-a_{0}|+|b_{1}-b_{0}|+|c_{1}-c_{0}|)
$$
for some constant $C>0$.
This inequality implies that if $a_{1}$, $b_{1}$, and $c_{1}$ approach $a_{0}$, $b_{0}$, and $c_{0}$ as $y\to 0$, then the Beltrami coefficient $b_{\widetilde{H}_{1}} (z)$ approaches to $b_{\widetilde{H}_{0}} (z)$ as $y\to 0$. Clearly if $\rho_{1}$ approaches $\rho_{0}$ as $y\to 0$, then $a_{1}$, $b_{1}$, and $c_{1}$ approach $a_{0}$, $b_{0}$, and $c_{0}$ as $y\to 0$. On the other hand, since 
$$
|b_{\widetilde{H}_{1}\circ \widetilde{H}_{0}^{-1}}(z)|\leq C |b_{\widetilde{H}_{1}}(z)-b_{\widetilde{H}_{0}}(z)| 
$$
for some constant $C>0$, if the Beltrami coefficient $b_{\widetilde{H}_{1}}(z)$ approaches the Beltrami coefficent $b_{\widetilde{H}_{0}}(z)$ as $y\to 0$, then the Beltrami coefficient $  
b_{\widetilde{H}_{1}\circ \widetilde{H}_{0}^{-1}}(z)$ approaches $0$ as $y\to 0$. 

Suppose $\Lambda (m)$ is the extremal length of the curve family that joins the interval $[-\infty, -1]$ to $[0,m]$, then it is an increasing real analytic function of $m$, in particular,
$\log m_{1}$ approaches to $\log m_{0}$ 
if and only if $\log \Lambda (m_{1})$ 
approaches $\log \Lambda( m_{0})$ 
(refer to~\cite[pages 74-76]{AhlforsBook}). 
The map $\widetilde{H}_{1}\circ \widetilde{H}_{0}^{-1}$ 
carries the extremal length problem for the family of curves joining $[-\infty, H_{0}(x-y)]$ to $[H_{0} (x), H_{0}(x+ky)]$ to the extremal length problem for the family of curves joining $[-\infty, H_{1}(x-y)]$ to $[H_{1} (x), H_{1}(x+ky)]$. If the Beltrami coefficient $b_{\widetilde{H}_{1}\circ \widetilde{H}_{0}^{-1}}(z)$ approaches $0$ as $y\to 0$, then the latter extremal length problem approaches the previous extremal length problem. This implies that 
$$
\frac{H_{1}(x+ky)-H_{1}(x)}{H_{1}(x)-H_{1}(x-y)}-
\frac{H_{0}(x+ky)-H_{0}(x)}{H_{0}(x)-H_{0}(x-y)}\to 0
$$
as $y\to 0$. From the above calculations, we proved a proposition that we will use in the proof of a complex characterization of uniform symmetry, as follows.

\medskip
\begin{proposition}~\label{skewqd}
Suppose $h_{0}$ and $h_{1}$ are two quasisymmetric circle homeomorphisms.
Then there is a bounded positive
real function $\epsilon (t)$ 
with $\epsilon (t)\to 0^{+}$ as $t\to 0^{+}$ such that
$$
|b_{\widetilde{H}_{1}} (z)-b_{\widetilde{H}_{0}} (z)| \leq \epsilon (y), \quad \forall z=x+iy \in \UH,
$$
if and only if there is another bounded positive real function $\widetilde{\epsilon}(t)$ with $\widetilde{\epsilon} (t)\to 0^{+}$ as $t\to 0^{+}$ such that 
$$
|\lambda_{H_{0}}(x,y,k) -\lambda_{H_{1}}(x,y,k)| \leq \widetilde{\epsilon} (|y|), \quad \forall x, y\in \R, \; 0\leq k\leq 1.
$$
\end{proposition}

Now we state our complex characterization result for a uniformly symmetric circle endomorphism. 

\medskip
\begin{theorem}~\label{ccusthm}
A circle endomorphism $f$ of degree $d$ is uniformly symmetric if and only if there are a quasiconformal extension $\widehat{F}$ of $F$ to the upper-half plane $\UH$ and a bounded positive real function $\epsilon (t)$ with $\epsilon (t)\to 0^{+}$ as $t\to 0^{+}$ such that 
$$
|b_{\widehat{F}^{-n}} (z)| \leq \epsilon (y), \quad \forall z=x+iy\in \UH, \quad \forall n>0.
$$
\end{theorem}

\begin{proof}
The ``if'' part is to apply the extremal length problem for the family of curves joining $[-\infty, x-y]$ to $[x, x+ky]$ to the extremal length problem for the family of curves joining $[-\infty, F^{-n} (x-y)]$ to $[F^{-n}(x), F^{-n} (x+ky)]$. Since 
$$
|b_{\widehat{F}^{-n}} (z)| \leq \epsilon (y), \quad \forall z=x+iy\in \UH, \quad \forall n>0,
$$
we have another bounded positive real function $\widehat{\epsilon}(t)$ with $\widehat{\epsilon}(t)\to 0^{+}$ as $t\to 0^{+}$ such that the quasisymmetric distortion function
$$
|\lambda_{F^{-n}} (x, y, k)-k|\leq \widehat{\epsilon} (|y|), \quad \forall x, y\in \UH, \; 0\leq k\leq 1.
$$
This implies that $F$ is uniformly symmetric. 

We now prove the ``only if'' part. 
Suppose $f$ is a uniformly symmetric circle endomorphism of degree $d$. Let $F$ be its lift. Then we have
$$
F(x+1) =F(x)+d, \quad \forall x\in \R.
$$
Let $h$ be the unique quasisymmetric homeomorphism with $h(1)=1$ such that
$$
f\circ h=h\circ q_{d}.
$$
Then we have 
$$
F\circ H =H\circ Q_{d}
$$
where $Q_{d} (z) =dz$ is the lift of $q_{d}(z)=z^{d}$.
Let $\widetilde{H}$ be the Beurling-Ahlfors extension of $H$ in our normalization. Then we have a quasiconformal extension $\widetilde{F}$ of $F$ to the upper-half plane $\UH$,
$$
\widetilde{F} =\widetilde{H}\circ Q_{d}\circ \widetilde{H}^{-1}.
$$
This implies that
$$
\widetilde{F}^{-n} =
\widetilde{H}\circ Q_{d}^{-n} \circ \widetilde{H}^{-1}.
$$
So we have that the Beltrami coefficient
\begin{equation}~\label{ceq}
|b_{\widetilde{F}^{-n}} (z)| =O(|b_{\widetilde{H}} (d^{-n}z) -b_{\widetilde{H}}(z)|)
\end{equation}
where $O(\cdot)$ is independent of $n$. 

Let $H_{0} =F^{-n}\circ H$ and $H_{1}=H$. Since $F$ is uniformly symmetric, we have a bounded positive real function $\widetilde{\epsilon} (t)$
with $\widetilde{\epsilon} (t)\to 0^{+}$ as $t\to 0^{+}$ such that 
quasisymmetric distortion functions satisfy
$$
|\lambda_{H_{1}}(x, y,k)- \lambda_{H_{0}} (x, y, k)| \leq \widetilde{\epsilon}(|y|), \quad \forall x, y\in \R, \; 0\leq k\leq 1. 
$$
Applying Proposition~\ref{skewqd},
there is another bounded 
positive function $\widehat{\epsilon} (t)$ 
with $\widehat{\epsilon} (t)\to 0^{+}$ 
as $t\to 0^{+}$ such that 
Beltrami coefficients satisfy
$$
|b_{\widetilde{F^{-n}\circ H}} (z) -b_{\widetilde{H}}(z)|\leq \widehat{\epsilon} (y), \quad \forall z=x+iy\in \UH,\; n\geq 1.
$$
Since 
$\widetilde{F^{-n}\circ H}= \widetilde{H\circ Q_{d}^{-n}}= \widetilde{H}\circ Q_{d}^{-n}$,
$$
b_{\widetilde{F^{-n}\circ H}} =b_{\widetilde{H}} (d^{-n} z)
$$ 
and 
$$
|b_{\widetilde{F}^{-n}} (z)| \leq C |b_{\widetilde{H}} (d^{-n}z) -b_{\widetilde{H}}(z)|
= C|b_{\widetilde{F^{-n}\circ H}}- b_{\widetilde{H}}(z)|
\leq C\widehat{\epsilon} (y)
$$
for all $z=x+iy\in \UH$ and all $n>0$.
Take $\epsilon (t) =C\widehat{\epsilon}(t)$. Then we have a bounded positive real function $\epsilon(t)$ with $\epsilon (t)\to 0^{+}$ as $t\to 0^{+}$ such that 
$$
|b_{\widetilde{F}^{-n}} (z)| \leq \epsilon (y), \quad \forall z=x+iy\in \UH,\; n\geq 1.
$$
We completed the proof.  
\end{proof}

\section{\bf Teichm\"uller Complex Manifold Structure}~\label{tcms}

Let $\mathcal{GM}$ be the space of all geometric Gibbs measures $\mu_{g}$ associated with $g\in \mathcal{G}$. From Theorem~\ref{main}, there is a bijective map 
$$
\eta: \mathcal{G} \to \mathcal{GM};\quad \eta (g)=\mu_{g}.
$$ 
As a subspace of $\mathcal{C}$, we have a maximum metric on $\mathcal{G}$, which induces a maximum metric on $\mathcal{GM}$,
$$
d_{max}(\mu_{1}, \mu_{2}) =\max_{w\in \Si} |\eta^{-1} (\mu_{1})(w) -\eta^{-1} (\mu_{2})(w)|, \quad \mu_{1}, \mu_{2}\in \mathcal{GM}
$$
However, this maximum norm does not capture the dynamical behavior of $(\Sigma, \s, g, \mu)$ well globally, especially when the geometric Gibbs measure is near the boundary. The induced Teichm\"uller metric 
$$
d_{\mathcal{GM}} (\mu_{1}, \mu_{2})= d_{\mathcal{G}} (\eta^{-1}(\mu_{1}), \eta^{-1}(\mu_{2})), \quad \mu_{1}, \mu_{2}\in \mathcal{GM},
$$
will measure the dynamical behavior of $(\Sigma, \s, g, \mu)$ well as a global hyperbolic metric. From Theorem~\ref{completion} and Theorem~\ref{main}, we get 

\medskip
\begin{corollary}~\label{completion2}
The Teichm\"uller space $(\mathcal{GM}, d_{\mathcal{GM}})$ is a complete space and the completion of $(\mathcal{HGM}, d_{\mathcal{HGM}})$, that is
$$
(\mathcal{GM}, d_{\mathcal{GM}})=\overline{(\mathcal{HGM}, d_{\mathcal{HGM}})}.
$$
\end{corollary}

As a consequence of Theorem~\ref{main} and Theorem~\ref{unifcont} 
\medskip
\begin{corollary}
The identity map
$$
Id_{TM}: (\mathcal{GM}, d_{\mathcal{GM}}(\cdot, \cdot)) \to (\mathcal{GM}, d_{max}(\cdot, \cdot))
$$
is uniformly continuous and the identity map
$$
Id_{MT}: (\mathcal{GM}, d_{max}(\cdot, \cdot))\to (\mathcal{GM}, d_{\mathcal{GM}}(\cdot, \cdot))
$$
is just continuous.
Thus the topology on $\mathcal{GM}$ induced by $d_{max} (\cdot, \cdot)$ 
and the topology on $\mathcal{GM}$ induced by $d_{\mathcal{GM}}(\cdot, \cdot)$ are the same. 
\end{corollary}

Let ${\mathcal H}$ be the space of all holomorphic periodic functions $\phi$ of period $1$ defined on the lower half plane ${\mathbb L}=\{ z=x+yi \in {\mathbb C}\;|\; y<0\}$ such that the norm 
$$
\|\phi \| = sup_{z\in {\mathbb L}} |y^{2} \phi (z)| <\infty.
$$
Let ${\mathcal H}_{0}$ be the subspace of ${\mathcal H}$ such that every element $\phi$ in it satisfies 
$$
|y^{2} \phi (z)| \to 0, \quad \hbox{as} \quad y\to 0 \quad \hbox{uniformly on $x$}.
$$
 
 For any quasisymmetric circle homeomorphism $h$, let $\widehat{h}$ be a quasiconformal extension of $h$ to the unit disk $\D$ and let $\widehat{H}$ be its lift on the upper-half plane ${\mathbb H}$. 
Let 
$$
b_{\widehat{H}}(z)  =\frac{\widehat{H}_{\overline{z}}(z)}{\widehat{H}_{z}(z)}
$$
be its Beltrami coefficient defined on ${\mathbb H}$. 
It is a periodic function of period $1$. Define
$$
b(z) =\left\{\begin{array}{ll}
b_{\widehat{H}} (z), &z\in {\mathbb H};\cr
0, &z\in {\mathbb L}
\end{array}\right.
$$
Let $w^{b}$ be the unique quasiconformal homeomorphism with Beltrami coefficient $b(z)$ that fixes $0$, $1$, and $\infty$; we call it a normalized quasiconformal homeomorphism. Then, it is analytic 
in ${\mathbb L}$ and $w^{b}(z+1)=w^{b}(z)+1$.
The Schwarzian derivative on $\mathbb{L}$ is 
$$
\phi (z) =S(w^{b}) (z) =\frac{(w^{b}) '''}{(w^{b})'} -\frac{3}{2} \Big( \frac{(w^{b})''}{(w^{b})'}\Big)^{2}.
$$
For two analytic functions $w^{b}$ and $w^{\widetilde{b}}$ on ${\mathbb L}$, we have a cocycle condition for Schwarzian derivatives, 
$$
S(w^{\widetilde{b}}\circ w^{b}) =S(w^{\widetilde{b}}) \circ w^{b}\cdot [(w^{b})']^{2} +S(w^{b}).
$$
Thus, if $w^{\widetilde{b}}$ is a M\"obius transformation, then we have $S(w^{\widetilde{b}}\circ w^{b})=S(w^{b})$. So $\phi$ is independent of the choice of any element in the equivalent class $[h]$ in the universal Teichm\"uller space $\mathcal{T}$. Bers' embedding 
$$
\iota: [h]\to \phi
$$  
embeds the universal Teichm\"uller space $\mathcal{T}$ into an open set in ${\mathcal H}$ containing the ball of radius $1/2$ centered at $0$ and contained in the ball of radius $3/2$ centered at $0$. 

Bers' embedding $\iota$ on the ball of radius $1/2$ centered at $0$ in ${\mathcal H}$ has a holomorphic section $s$ by Ahlfors-Weill's technique (see~\cite {AhlforsBook}). 
For any analytic function $\phi \in {\mathcal H}$ with $\|\phi\|<1/2$. Consider two independent solutions $\eta_{1}$ and $\eta_{2}$ of the ordinary differential equation
$$
\eta''=-\frac{1}{2} \phi  \eta.
$$ 
Define 
$$
 G (z) = \left\{ \begin{array}{ll}
\frac{\eta_{1}(\overline{z}) +(z-\overline{z}) \eta_{1}' (\overline{z})}{\eta_{2}(\overline{z}) +(z-\overline{z}) \eta_{2}' (\overline{z})}, & z\in {\mathbb H};\cr
\frac{\eta_{1}(z)}{\eta_{2}(z)}, & z\in {\mathbb L}.
\end{array}\right.
$$
The Beltrami coefficient of $G(z)$ is
$$
b_{G} (z) =\left\{ \begin{array}{ll} 
\frac{1}{2} (z-\overline{z})^{2} \phi (z), & z\in {\mathbb H};\cr
0, & z\in {\mathbb L}.
\end{array}\right.
$$
Thus $\|b_{G}\|_{\infty} <1$ and is a periodic function of period $1$.
Define 
$$
\widetilde{b} (z) =\left\{
\begin{array}{ll}
b_{G} (z), & z\in {\mathbb H};\cr
\overline{b_{G} (\overline{z})}, & z\in {\mathbb L}.
\end{array}\right.
$$
It is still a periodic function of period $1$. 
Let $w_{\widetilde{b}}$ be the normalized quasiconformal homeomorphism with the Beltrami coefficient $\widetilde{b}$. Then $H=w_{\widetilde{b}}|{\mathbb R}$ is a homeomorphism of ${\mathbb R}$ with $H(x+1) =H(x)+1$. So it is the lift of a circle homeomorphism $h$. 
The local section 
$$
s: \phi \to [h]
$$  
is the local right inverse of $\iota$ defined on a small disk $\Delta$ centered at $0$ in ${\mathcal H}$, that is, 
$\iota\circ s=id$ on $\Delta$.  Bers' embedding $\iota$, the local section $s$, and changing-basepoint isomorphisms give an atlas of holomorphic charts on $\mathcal{T}$, making $\mathcal{T}$ a complex manifold. 
 
According to~\cite{GS}, if $h$ is a symmetric circle homeomorphism, then Bers' embedding $\iota ([h])=\phi \in {\mathcal H}_{0}$ and the holomorphic section $s$ maps the ball of radius $1/2$ centered at $0$ in ${\mathcal H}_{0}$ into the space ${\mathcal S}$ of all symmetric circle homeomorphisms modulo M\"obius transformations. Thus Bers' embedding induces an embedding, which we still denote as $\iota$, from the universal asymptotical Teichm\"uller space $\mathcal{AT}=\mathcal{T}/\mathcal{S}$ into the quotient space ${\mathcal H}/{\mathcal H}_{0}$. The local section $s$ induces a local section, which we still denote as $s$ defined on the quotient space $\Delta/{\mathcal H}_{0}$ such that $\iota \circ s=id$. 
The embedding $\iota$, the local section $s$, and changing basepoint isomorphisms give holomorphic charts on $\mathcal{AT}$. Thus $\mathcal{AT}$ is also a complex manifold.  
What we are going to prove is that Theorem~\ref{main} allows us to follow a similar idea to Bers' embedding to obtain the Teichm\"uller complex manifold structure 
on $ \mathcal{GM}$ by using open sets in the quotient space ${\mathcal H}_{1}/{\mathcal H}_{0}$, 
where 
$$
{\mathcal H}_{1}=\{ \phi \in {\mathcal H}\;|\; y^{2} (\phi (d^{-n}z) -\phi (z))  \to 0 \; \hbox{as $y\to 0$ uniformly on $x$ and $n$}\},
$$
where $z=x+iy$.

From Theorem~\ref{main}, we have that 
$$
\beta=\eta\circ \zeta: \mathcal{TUS}\to \mathcal{GM}
$$ 
is a bijective map. For every $\mu\in \mathcal{GM}$, let 
$\tau=[(f, h_{f})]=\beta^{-1}(\mu)\in \mathcal{TUS}$ and 
$$ 
F=H\circ Q_{d}\circ H^{-1},
$$
where $Q_{d} =dz$ and $H=H_{f}$ are lifts of $q_{d}$ and $h_{f}$ to $\mathbb{R}$. Suppose $\widetilde{H}$ is the Beurling-Ahlfors extension of $H$ and 
and  
$$
\widetilde{F} =\widetilde{H}\circ Q_{d} \circ \widetilde{H}^{-1}.
$$
Let 
$$
b_{\widetilde{H}}(z)=\frac{\widetilde{H}_{\overline{z}}(z)}{\widetilde{H}_{z}(z)}
$$
be the Beltrami coefficient of $\widetilde{H}$ and let
$$
b_{\widetilde{F}^{-n}} (z) =\frac{\widetilde{F}^{-n}_{\overline{z}}(z)}{\widetilde{F}^{-n}_{z}(z)}
$$
be the Beltrami coefficient of $\widetilde{F}^{-n}$. 
From Theorem~\ref{ccusthm}, we have that 
$$
b_{\widetilde{H}}(d^{-n} z)-b_{\widetilde{H}}(z) \to 0 \quad \hbox{as $y\to 0$ uniformly on $x$ and $n$}
$$
and 
$$
b_{\widetilde{F}^{-n}}(z) \to 0 \quad \hbox{as $y\to 0$ uniformly on $x$ and $n$}.
$$

Define 
$$
\widehat{b} (z) =\left\{ \begin{array}{ll} 
b_{\widetilde{H}} (z), & z\in {\mathbb H};\cr
0, & z\in {\mathbb L}.
\end{array}\right.
$$
Let $w^{\widehat{b}}$ be the normalized quasiconformal map of the Riemann sphere $\CC$ whose Beltrami coefficient is $\widehat{b}$. Then $w^{\widehat{b}}$ is analytic on the lower-half plane $\mathbb{L}$. Let $\phi  =S(w^{\widehat{b}})\in {\mathcal H}$ be Bers' embedding of $h_{f}$. For any other representations of $[(f, h_{f})]$, the difference between Bers' embeddings of them and $\phi$ are in ${\mathcal H}_{0}$.
So we actually embed the equivalence class $[(f, h_{f})]$ into an equivalence class $[\phi]\in {\mathcal H}/{\mathcal H}_{0}$. We still use $\iota$ to denote this embedding. 
We further define 
$$
b_{n} (z)=\left\{ \begin{array}{ll} 
b_{\widetilde{F}^{-n}} (z), & z\in {\mathbb H};\cr
0, & z\in {\mathbb L}.
\end{array}\right\}
$$
Let $w^{b_{n}}$ be the normalized quasiconformal map of $\CC$ whose Beltrami coefficient is $b_{n}$.
Since 
$
w^{b_{n}}=w^{\widehat{b}}\circ Q_{d}^{-n}\circ (w^{\widehat{b}})^{-1},
$
$$
S(w^{b_{n}})((w^{\widehat{b}})^{-1}(z)) = S(w^{\widehat{b}})(d^{-n}z) -S (w^{\widehat{b}})(z) =\phi(d^{-n}z)-\phi (z). 
$$
Since $b_{\widehat{F}^{-n}}(z) \to 0$ as $y\to 0$ uniformly on $x$ and $n$, we get $\phi \in {\mathcal H}_{1}$ and, furthermore,  we have that 
$$
\iota (h_{f}) = [\phi]\in {\mathcal H}_{1}/{\mathcal H}_{0}.
$$
 
For any $[\phi] \in {\mathcal H}_{1}/{\mathcal H}_{0}$ with $\|\phi\|<1/2$ and the corresponding $h=s(\phi)$, 
$$
f=h\circ q_{d}\circ h^{-1}
$$ 
is a uniformly symmetric circle endomorphism. Let $\mu= \beta ([f, h])$. It is independent of the choice of representation in $[\phi]$. Thus the local section $s$ induces a local section, which we still denote as $s$ defined on the open ball $\widehat{U}$ of radius $1/2$ centered at $0$ in ${\mathcal H}_{1}/{\mathcal H}_{0}$ such that 
 $\iota \circ s=id$.
 The embedding $\iota$, the local section $s$, and the changing-basepoint isomorphisms give an atlas of holomorphic charts on $\mathcal{GM}$. Thus $\mathcal{GM}$ is a complex manifold, and we have proved that  

\medskip
\begin{starthmb}~\label{cms}
The space $(\mathcal{GM}, d_{\mathcal{GM}})$ is a complete complex manifold modeled
on ${\mathcal H}_{1}/{\mathcal H}_{0}$.
\end{starthmb}

Remember that $\mu_{0}$ is the basepoint of $\mathcal{GM}$. The tangent space $T_{\mu_{0}}\mathcal{GM}$ of $\mathcal{GM}$ at the basepoint can now be thought of as a subspace of the tangent space $T_{[id]}\mathcal{AT}$ of the universal asymptotic Teichm\"uller space $\mathcal{AT}$ at the basepoint as follows.  

Consider a smooth path $[h_{t}]$ in the universal Teichm\"uller space $\mathcal{T}$ such that $h_{0}=id$. Let $H_{t}$ be the corresponding path of lifts. 
According to~\cite{Ri}, the tangent vector 
$$
E (x) =\frac{dH_{t}(x)}{dt}|_{t=0}
$$
is a periodic Zgymund continuous function of period $1$ vanishing at all integers $n$, that is,
\begin{equation}~\label{bz}
E(x+y) +E(x-y)-2E(x) = O(y), \quad E(x+1)=E(x), \quad E(n)=0,
\end{equation}
Let $\mathcal{Z}$ be the space of all Zygmund continuous functions $E$ on the real line $\R$ satisfying (\ref{bz}). 
Consider a smooth path $[h_{t}]$ in $\mathcal{S}$ such that $h_{0}=id$, then $E$ is a little  
Zgymund continuous function of period $1$ vanishing at all integers $n$, that is,
\begin{equation}~\label{lz}
E(x+y) +E(x-y)-2E(x) = o(y), \quad E(x+1)=E(x), \quad E(n)=0.
\end{equation}
Let $\mathcal{Z}_{0}$ be the subspace of $\mathcal{Z}$ consisting of all little Zygmund continuous functions on the real line $\R$ satisfying (\ref{lz}). 
Then 
$$
T_{[id]}\mathcal{AT} = \mathcal{Z}/\mathcal{Z}_{0}.
$$
 
Now consider a path $\mu_{t}$ in $\mathcal{GM}$ passing $\mu_{0}$ at $t=0$. Let $[(f_{t}, h_{t})]=\beta^{-1}(\mu_{t}) \in \mathcal{TUS}$ passing the basepoint $[(q_{d}, id)]$. 
Let $F_{t}$ and $H_{t}$ be the lift of $f_{t}$ and $h_{t}$. Then we have $E\in \mathcal{Z}$ such that  
$$
H_{t} (x) = x +t E(x) +o(t)\quad \hbox{and}\quad H^{-1}_{t}(x) =x-t E(x) +o(t).
$$
Now consider 
$$
F_{t}^{-n} (x) = H_{t}\circ Q_{d}^{-n} \circ H_{t}^{-1} (x)= H_{t} (d^{-n} x -td^{-n} E(x)+o(t))
$$
$$
= d^{-n} x-td^{-n} E(x) +o(t) +t E (d^{-n}x -td^{-n} E(x) +o(t))
$$
$$
= d^{-n} x + (E(d^{-n}x) -d^{-n} E(x)) t +o(t)=d^{-n}x +V_{n} (x) t +o(t) 
$$
where 
\begin{equation}~\label{ulz}
V_{n} (x) = E(d^{-n}x) -d^{-n} E(x).
\end{equation}
Since all $F_{t}^{-n}$ are uniformly symmetric, according to~\cite{Ri}, functions $V_{n}(x)$ are uniformly little Zygmund continuous functions, that is, there is a bounded positive continuous real function $\varepsilon (t)$ with $\varepsilon (t) \to 0^{+}$ as $t\to 0^{+}$ such that
\begin{equation}~\label{ulz1}
|V_{n}(x+y) +V_{n}(x-y) - 2V_{n}(x) |  < \varepsilon (y)y
\end{equation}
for all $x\in \mathbb{R}$, $y>0$, and all $n$, and moreover, 
\begin{equation}~\label{ulz2}
V_{n}(x+d^{n}) =V_{n} (x).
\end{equation}
Let $\mathcal{Z}_{1}$ be the space of all Zgymund continuous functions $E$ on the real line $\R$ satisfying (\ref{bz}), (\ref{ulz1}), and (\ref{ulz2}). We have proved that 

\medskip
\begin{starthmc}~\label{tg} 
$$
T_{\mu_{0}} \mathcal{GM} = \mathcal{Z}_{1}/\mathcal{Z}_{0}.
$$
\end{starthmc}
    
\medskip
The Hilbert transformation ${\mathcal J}$ on ${\mathcal Z}$ is defined by
$$
{\mathcal J}E(x)=\frac{1}{\pi}\int^{\infty}_{-\infty}E(y)R(x,y) dy
$$
where
$$
R(x,y)=\frac{1}{y-x}-\frac{x}{y-1}+\frac{x-1}{y}.
$$
Let $\zeta=\xi+i\eta$ be a complex number. Then $\overline{\zeta}=\xi-i\eta$. By Stokes' formula,
$$
{\mathcal
J}E(x)=\frac{1}{\pi}\int^{\infty}_{-\infty}E(y)R(x,y)
dt=\frac{2i}{\pi}\int\int_{\mathbb H}
\overline{\partial}\widehat{E}(\zeta) R(x,\zeta)d\xi d\eta+iE(x),
$$
where ${\mathbb H}$ is the upper-half plane and $\widehat{E}$ is a continuous extension of $E$ to the complex plane and 
$\overline{\partial}\widehat{E} = {\partial \widehat{E}}/{\partial\overline{\zeta}}$.

Let $\mu (\zeta)=\overline{\partial}\widehat{E}(\zeta)$, then
$$
{\mathcal J}E(x)=\frac{2i}{\pi}\int\int_{\mathbb H}
\mu(\zeta) R(x,\zeta)d\xi d\eta+iE(x)
$$
$$
=\frac{-2i}{\pi}\int\int_{\mathbb L}
\mu(\zeta)R(x,\zeta)d\xi d\eta - iV(x)
$$
where ${\mathbb L}$ is the lower-half plane.

Define
$$
\widetilde{\mu} (\zeta) =\left\{
                              \begin{array}{rr}
                                    -i\mu (\zeta), & \zeta\in {\mathbb H};\\
                                     i\mu (\zeta), & \zeta\in {\mathbb L}.
                              \end{array}\right.
$$
Then
$$
{\mathcal J}E(x)=\frac{-1}{\pi}\int\int_{\mathbb C}
\widetilde{\mu}(\zeta)R(x,\zeta)d\xi d\eta.
$$
And we also have
$$
\overline{\partial}({\mathcal J}E)=-i\mu\quad \hbox{and}\quad
\overline{\partial}({\mathcal J}^2E)=i^2\mu=-\mu.
$$
This implies
$$
{\mathcal J}^2E=-E.
$$
Hence, the Hilbert transformation ${\mathcal J}$ gives an almost complex structure on ${\mathcal Z}$.

Now we will prove that the Hilbert transformation ${\mathcal J}$ also gives an almost
complex structure on the tangent space $T_{\mu_{0}} \mathcal{GM}$. 

\vspace*{5pt}
\begin{starthmd}~\label{acs}
Let $[{\mathcal J}E]$ be the Hilbert transformation of $[E] \in {\mathcal
Z}_{1}/ {\mathcal Z}_{0}$. Then ${\mathcal J}E$ also satisfies (\ref{bz}), (\ref{ulz1}), and (\ref{ulz2}). 
Therefore, ${\mathcal J}[E]:=[{\mathcal J}E] \in {\mathcal
Z}_{1}/ {\mathcal Z}_{0}$. Since
${\mathcal J}: {\mathcal
Z}_{1}/ {\mathcal Z}_{0}\to {\mathcal
Z}_{1}/ {\mathcal Z}_{0}$ and ${\mathcal
J}^{2}=-{\mathcal I}$, it is an almost complex structure on
$T_{\mu_{0}}\mathcal{GM}$.
\end{starthmd}

\begin{proof}
Note that
$$
{\mathcal J}E(x)=\frac{1}{\pi}\int^{\infty}_{-\infty}E(y)\Big( \frac{1}{y-x}-\frac{x}{y-1}+\frac{x-1}{y}\Big) dy.
$$
From $E(x+1)=E(x)$, we have
$$
\frac{1}{\pi}\int^{\infty}_{-\infty}\frac{x E(y)}{y-1} dy
=\frac{1}{\pi}\int^{\infty}_{-\infty}\frac{xE(y+1)}{y}
dy =\frac{1}{\pi}\int^{\infty}_{-\infty}\frac{xE(y)}{y}
dy.
$$
Hence
$$
{\mathcal J}E(x)=\frac{1}{\pi}\int_{-\infty}^{+\infty}\frac{E(y)}{y-x} dy-
\frac{1}{\pi}\int^{\infty}_{-\infty}\frac{E(y)}{y}
dy.
$$
From this form of the Hilbert transform, it is easy to check
${\mathcal J}E(0)=0$ and ${\mathcal J}E(x+1)={\mathcal J}E(x)$ since $E(y+1)=E(y)$.
Now
$$
\widetilde{V}_{n}(x) ={\mathcal J}E (d^{-n}x) -d^{-n} {\mathcal J} E (x)
$$
$$
=\frac{1}{\pi}\int_{-\infty}^{+\infty}\frac{E(y)}{y-d^{-n}x} dy-
\frac{1}{\pi}\int_{-\infty}^{+\infty}\frac{d^{-n} E(y)}{y-x} dy-
(1-d^{-n}) \frac{1}{\pi}\int^{\infty}_{-\infty}\frac{E(y)}{y}
$$
$$
=\frac{1}{\pi}\int_{-\infty}^{+\infty}\frac{E(d^{-n}y)-d^{-n} E(y)}{y-x} dy-
(1-d^{-n}) \frac{1}{\pi}\int^{\infty}_{-\infty}\frac{E(y)}{y}
$$
$$
=\frac{1}{\pi}\int_{-\infty}^{+\infty}\frac{V_{n}(y)}{y-x} dy-
(1-d^{-n}) \frac{1}{\pi}\int^{\infty}_{-\infty}\frac{E(y)}{y}.
$$
The last formula implies that $\widetilde{V}_{n}$ is uniformly little Zygmund continuous since $V_{n}$ is uniformly little Zygmund continuous. 
We proved the theorem.
 \end{proof}

\section{Welding}~\label{weld}

Suppose $\mu\in \mathcal{GM}$ is a geometric Gibbs measure and $g\in \mathcal{G}$ is the corresponding continuous potential. 
Thus we have an $f_{\mu}\in \mathcal{TUS}$ such that $\tau (f_{\mu})=g$.
The corresponding 
$$
\mu \longleftrightarrow g\longleftrightarrow f_{\mu}
$$
is bijective. Furthermore, we have a unique quasisymmetric circle homomorphism $h_{\mu}$ with $h_{\mu}(1)=1$ such that
$$
f_{\mu} =h_{\mu}\circ q_{d}\circ h_{\mu}^{-1}.
$$
Using $h=h_{\mu}$, we can weld $f=f_{\mu}$ and $q_{d}$ to get a dynamical system $\alpha$ on the Riemann sphere $\CC$ such that it preserves a quasicircle, which we will call a dynamical uniformly asymptotically conformal quasicircle or just a dynamical quasicircle.  

Consider the Beurling-Ahlfors extension $\widetilde{h}$ of $h$ to the unit disk $\D$ and its Beltrami coefficient 
$$
b_{\widetilde{h}} (z) =\frac{\widetilde{h}_{\overline{z}}(z)}{\widetilde{h}_{z}(z)}, \quad z\in \D .
$$
Define 
$$
\widehat{b} (z) =\left\{ \begin{array}{ll} 
b_{\widetilde{h}} (z), & z\in \D;\cr
0, & z\in \CC\setminus \D.
\end{array}\right.
$$
Let $w^{\widehat{b}}$ be the normalized quasiconformal mapping of $\CC$ whose Beltrami coefficient is $\widehat{b}$. 
Then 
$$
D=D_{\mu}=w^{\widehat{b}}(\D)
$$
is a quasidisk and 
$$
S=S_{\mu}= w^{\widehat{b}}(\T)=\partial w^{\widehat{b}}(\D)
$$
is a quasicircle, which we call a dynamical quasicircle corresponding to $\mu$. 

The map 
$$
R_{\infty} =(w^{\widehat{b}})^{-1}: 
\CC\setminus \overline{D}\to \CC\setminus \overline{\D}
$$
is the Riemann mapping with $R_{\infty} (\infty)=\infty$ and $R_{\infty}(1)=1$. Let 
$$
R_{0}: D\to \D
$$ 
be the Riemann mapping with $R_{0}(0)=0$ and $R_{0}(1)=1$. 
Both $R_{0}$ and $R_{\infty}$ can be extended to $\overline{D}$ and $\CC\setminus D$ as homeomorphisms. Then 
$$
R_{0}\circ R_{\infty}^{-1}: \T\to \T
$$ 
is a homeomorphism of the unit circle $\T$ fixing $1$. 
Since $\tilde{h}^{-1}\circ R_{0}$ and $(w^{\widehat{b}})^{-1}$ have the same Beltrami coefficient on $D$ and both are normalized, we have that 
$$
\tilde{h}^{-1} \circ R_{0} =(w^{\widehat{b}})^{-1}
$$
on $D$ and
$$
h^{-1}\circ R_{0} =R_{\infty}
$$
on $\partial D$. Thus, we have that
$$
h= R_{0} \circ R_{\infty}^{-1}
$$
on $\T$.

Since $f$ is uniformly symmetric, from Theorem~\ref{ccusthm}, we have a quasiconformal extension $\widehat{f}$ of $f$ to $\CC$ such that 
$\widehat{f}$ has $0$ as an attracting fixed point and such that
$$
b_{\widehat{F}^{-n}} (z) \to 0, \quad \hbox{as} \quad y\to 0 \quad \hbox{uniformly on $x\in \R$ and $n\geq 1$}
$$
where $\widehat{F}$ is the lift of $\widehat{f}$. 
Define
$$
\alpha (z) =\left\{ \begin{array}{ll} 
R_{0}^{-1}\circ \widehat{f}\circ R_{0} (z), & z\in D;\cr
R_{\infty}^{-1} \circ q_{d} \circ R_{\infty}, & z\in \CC\setminus D.
\end{array}\right.
$$
Then $\alpha$ is a continuous branched cover of $\CC$ with only two branched points $0$ and $\infty$ such that $D$ is the basin of $0$ and $\CC\setminus \overline{D}$ is the basin of $\infty$, that is,
$$
D=\{z\in \CC\;|\; \alpha^{n} (z) \to 0\hbox{ as $n\to \infty$}\}
$$
and 
$$
\CC\setminus \overline{D}=\{z\in \CC\;|\; \alpha^{n} (z) \to \infty \hbox{ as $n\to \infty$}\}
$$
The quasicircle $S$ is completely invariant under $\alpha$, that is, 
$$
\alpha (S)=\alpha^{-1}(S)=S.
$$
The branched cover $\alpha$ of degree $d$ uniformizes $\widehat{f}$ and $q_{d}$ simultaneously, that is,
$$
\widehat{f} = R_{0} \circ \alpha\circ R_{0}^{-1}
$$
on $\overline{D}$ and 
$$
q_{d} =R_{\infty} \circ \alpha\circ R_{\infty}^{-1}
$$
on $\CC\setminus D$. 

Recall the symbolic space $\Sigma$ in Section 1. For any $w=\cdots i_{n-1}\cdots i_{1}i_{0}\in \Si$, we have $w_{n}=i_{n-1}\cdots i_{1}i_{0}$, $\sigma(w_{n})=i_{n-1}\cdots i_{1}$, $J_{w_{n}}$, and $J_{\sigma(w_{n})}$ in Section 1. 
Let $S_{w_{n}}=R_{0}^{-1} (J_{w_{n}})$ and $S_{\sigma(w_{n})}=R_{0}^{-1} (J_{\sigma(w_{n})})$. 

\medskip
\begin{starthme}~\label{hmgm}
Suppose $\nu$ is the harmonic measure for the marked quasidisk $(D, 0)$. Then it is an $\alpha$-invariant measure, that is, $\alpha_{*}\nu=\nu$, and a geometric Gibbs measure in the meaning that
$$
\lim_{n\to \infty} \frac{ \nu (S_{w_{n}})}{\nu (S_{\sigma(w_{n})})} =g(w)
$$
where $g$ is a continuous function on $\Sigma$. The convergence is uniform on $\Si$. When the convergence is exponential, then $g$ is a H\"older continuous function on $\Sigma$ and $\nu$ is the Gibbs measure associated with $g$ in the classical sense. 
\end{starthme}

\begin{proof}
Let $l(J)=|J|$ be the Lebesgue measure on $\T$. The harmonic measure $\nu=(R_{0}^{-1})_{*} l$, that is, $\nu (A) =l(R_{0} (A))$ for any measurable set $A\subset S$. Since $f$ preserves the Lebesgue measure $l$ and $\alpha = R_{0}^{-1}\circ f\circ R_{0}$ on $S$, $\nu$ is $\alpha$-invariant, and for $w\in \Sigma$,
$$
\lim_{n\to \infty} \frac{ \nu (S_{w_{n}})}{\nu (S_{\sigma(w_{n})})} =
\lim_{n\to \infty} \frac{|J_{w_{n}}|}{|J_{\sigma(w_{n})}|} 
\lim_{n\to \infty} \frac{\mu([w]_{n})}{\mu ([\sigma ( w)]_{n-1})} 
=g(w).
$$
The convergence is uniform on $\Si$. 

When the convergence is exponential, $g$ is a H\"older continuous function potential, and we can find a constant $C\geq 1$ such that 
$$
C^{-1} \leq \frac{\nu (S_{w_{n}})}{\prod_{k=0}^{n-1} g(\sigma^{k}(w))}\leq C, \quad \forall w\in \Si, \;n\geq 1.
$$
Thus, it is the Gibbs measure in the classical sense.
\end{proof}

Suppose $\mu_{j}\in \mathcal{GM}$, for $j=0,1$, and 
suppose $S_{j}=S_{\mu_{j}}$ is the corresponding dynamical quasicircle. 
Let $\widehat{f}_{j}$ be an extension of $f_{\mu_{j}}$ to $\CC$ with 
$$
b_{\widehat{F}_{j}^{-n}} (z) \to 0, \quad \hbox{as} \quad y\to 0 \quad \hbox{uniformly on $x\in \R$ and $n\geq 1$}
$$
where $\widehat{F}_{j}$ is the lift of $\widehat{f}_{j}$. 
Let $\alpha_{j}$ be the branched cover of degree $d$ uniformizing $\widehat{f}_{j}$ and $q_{d}$ simultaneously. Then $\alpha_{j} (S_{j})=S_{j}$ and
$$
R_{0,j}\circ \alpha_{j}\circ R_{0, j}^{-1
} =\widehat{f}_{j}
$$
on $\D$ and
$$
R_{\infty,j}\circ \alpha_{j}\circ R_{\infty, j}^{-1
} =q_{d}
$$
on $\CC\setminus \D$.

We say $S_{0}$ and $S_{1}$ are asymptotically conformal equivalent if there is a normalized quasiconformal mapping $\widehat{h}$ of $\CC$ such that $\widehat{h} (S_{0})=S_{1}$ with the Beltrami coefficient $b_{\widehat{h}}(z) \to 0$ as $z\to S_{0}$ and such that 
$$
\alpha_{1} =\widehat{h}\circ \alpha_{0}\circ \widehat{h}^{-1}
$$
on $S_{1}$. The next theorem says that every dynamical quasicircle is asymptotically conformal rigid. 

\medskip
\begin{starthmf}~\label{ace}
If two dynamical quasicircles $S_{0}$ and $S_{1}$ are asymptotically conformal equivalent, then they are the same. 
\end{starthmf} 
 
\begin{proof}
Let $R_{0,j}: D_{j}\to \D$ be the Riemann mapping for $j=0,1$ fixing $0$ and $1$.
Define 
$$
h= R_{0, 1}\circ \widehat{h}\circ R_{0,0}^{-1}: \T\to \T.
$$
Then we have $h(1)=1$. On $\T$, we have that 
$$
f_{\mu_{1}} = R_{0, 1} \circ \alpha_{1}\circ R_{0,1}^{-1}
= R_{0, 1} \circ \widehat{h}\circ \alpha_{0}\circ \widehat{h}^{-1}\circ R_{0,1}^{-1} =h\circ f_{\mu_{0}} \circ h^{-1}.
$$
Since $R_{0, 1}\circ \widehat{h}\circ R_{0,1}^{-1}: \D\to \D$ is a quasiconformal map with the Beltrami coefficient
$$
b_{R_{0, 1}\circ \widehat{h}\circ R_{0,0}^{-1}}(z) \to 0 \quad \hbox{ as $z\to \T$},
$$
$h$ is a symmetric circle homeomorphism of $\T$. Since both $f_{\mu_{0}}$ and $f_{\mu_{1}}$ preserve the Lebesgue measure, Theorem A implies that $h=Id$. Thus, $S_{0}=S_{1}$. We proved the theorem.

\end{proof}

Let $\mathcal{DQC}$ be the space of all dynamical (uniformly asymptotically conformal) quasicircles. Theorem F implies that $\mathcal{DQC}$ and $\mathcal{GM}$ are in bijection. Theorem B implies that $\mathcal{DQC}$ is a complex manifold. 

For every dynamical quasicircle $S_{\mu}$ with $\mu \in \mathcal{GM}$,  
there exist three geometric Gibbs measures associated with three different potentials.  

\medskip
\noindent
1. \textbf{Harmonic measure at the origin.}  
The first is the harmonic measure 
\[
\nu = (R_{0}^{-1})_{*} l
\] 
for the marked quasidisk $(D_{\mu}, 0)$ on $S_{\mu}$.  
By Theorem E, this measure is the geometric Gibbs measure associated with a potential $g$, viewed as the dual derivative of $\alpha: S_{\mu} \to S_{\mu}$ defined on $\Sigma$.  

\medskip
\noindent
2. \textbf{Hausdorff measure.}  
The second is the Gibbs measure associated with the potential 
\[
\phi(z) = - \log d \, |R_{\infty}^{-1}(z)|, \qquad z \in \mathbb{T}.
\]  
Since $R_{\infty}^{-1}$ extends to a quasiconformal homeomorphism of the Riemann sphere $\mathbb{C}$, the potential $\phi$ is H\"older continuous on $\mathbb{T}$.  
The system $(q_{d}, \phi)$ admits a unique Gibbs measure $\nu_{\infty}$, which is equivalent to the Hausdorff measure $H^{t}(\cdot)$ for $t = \dim_{H}(S_{\mu})$, where $t$ is the unique solution of the pressure equation $P(t\phi) = 0$.  

\medskip
\noindent
3. \textbf{Harmonic measure at infinity.}  
The third is the harmonic measure 
\[
\nu_{q} = (R_{\infty}^{-1})_{*} l
\] 
on $S_{\mu}$ for the marked quasidisk $(\mathbb{C} \setminus \overline{D}_{\mu}, \infty)$.  
Since $R_{\infty}$ also extends to a quasiconformal homeomorphism of $\mathbb{C}$, the function  
\[
\phi_{d}(z) = - \log d \, |R_{\infty}(z)|, \quad z\in S_{\mu},
\]  
is H\"older continuous.  
Thus the system $(\alpha, \phi_{d})$ has a unique Gibbs measure $\nu_{d}$ on $S_{\mu}$.  
Moreover, since $(R_{\infty})_{*}\nu_{d} = l$ and $l$ is the Gibbs measure for the system $(q_{d}|\mathbb{T}, -\log d)$, it follows that  
\[
h_{\nu_{d}}(\alpha|S_{\mu}) = h_{l}(q_{d}|\mathbb{T}) = \log d.
\]  
Therefore, $\nu_{d}$ is the measure of maximal entropy (mme) for $\alpha|S_{\mu}$.

\section{\bf Ergodicity}~\label{erg} 

The ergodicity properties of geometric Gibbs measures remain an open area. In this section, we present a foundational result in this direction. 

\medskip
\begin{definition}~\label{uac}
Suppose $f$ is a circle endomorphism. We say that $f$ is uniformly absolutely continuous 
if all inverses branches of $f^{m}: J \to I$ are uniformly absolutely continuous with respect to conditional probability measures 
$\mu_{0} (\cdot | I)$ and $\mu_{0} (\cdot |J)$ for all intervals $J$ and $I$ in $\T$ and all integers $m\geq 1$ such that $f^{m}: J\to I$ are homeomorphisms, that is,  
for any $\epsilon >0$, there is a $\delta>0$ such that for any
countably many pairwise disjoint open intervals $\{ I_{i}\}$ in $I$ with 
$$
\frac{\sum_{i} |I_{i}|}{|I|} <\delta,
$$
we have 
$$
\frac{|f^{-m}(\cup_{i}I_{i}) \cap J|}{|J|} <\epsilon\quad \forall\; m\geq 1.
$$
\end{definition}

Let 
$$
\mathcal{UAC}=\{ f\in \mathcal{US}\;|\; f \hbox{ is uniformly absolutely continuous}\}.
$$
We can similarly define the Teichmüller space $\mathcal{TUAC}$ to that in Section 1. Then $\mathcal{TUAC}$ is a subspace of $\mathcal{TUS}$. An example of a uniformly absolutely continuous circle endomorphism is a $C^{1+Dini}$ expanding circle endomorphism. Suppose $f$ is a $C^{1}$ circle endomorphism. Then $f'$ is a continuous function and has a modulus of continuity 
$$
\omega (t) = \sup_{z,w\in \T, |z-w|\leq t} |f'(z)-f'(w)|.
$$
We say that $f$ is $C^{1+Dini}$ if 
$$
\int_{0}^{a} \frac{\omega(t)}{t} dt <\infty.
$$
If $\omega (t)\leq Ct^{\alpha}$ for some $0<\alpha\leq 1$, we say $f$ is $C^{1+\alpha}$. 
One can check that a $C^{1+\alpha}$ circle endomorphism is a $C^{1+Dini}$ circle endomorphism.   
 
\medskip
\begin{example}
Suppose $f$ is a $C^{1+Dini}$ circle expanding circle endomorphism. Then $f$ is in $\mathcal{US}$ and is uniformly absolutely continuous.
\end{example} 

\begin{proof}
Let 
$$
\widetilde{\omega}(t) =\sum_{n=1}^{\infty} \omega
(\lambda^{-n}t).
$$
Since $\omega(t)$ is Dini,  $\widetilde{\omega}(t) \to 0$ as $t\to 0$. 

Let $F$ be the lift of $f$ into the real line $\R$. We have that $F(x+1) = F(x)+d$, $F'(x+1)=F'(x)$, and 
$(F^{n})' (x) \geq C \lambda^{n}$ for all $n>0$
where $C>0$ and $\lambda>1$ are two constants for $f$. Without loss of generality, we assume that $C=1$. 

For any $x, y\in {\mathbb R}$ and $n>0$, let $x_{k}=F^{-k}(x)$ and
$y_{k}=F^{-k}(y)$, $0\leq k\leq n$. Then
$$
\Big| \log \frac{(F^{-n})'(x)}{(F^{-n})'(y)} \Big|= \Big| \log
\frac{(F^{n})'(y_{n})}{(F^{n})'(x_{n})}\Big| \leq \sum_{k=1}^{n} |
\log F'(x_{k})-\log F'(y_{k})|
$$
$$
\leq \frac{1}{\lambda} \sum_{k=1}^{n} | F'(x_{k})-F'(y_{k})|
\leq \frac{1}{\lambda} \sum_{k=1}^{n} \omega (
\lambda^{-k} |x-y|).
$$
Let $D=1/\lambda$. Then we have the following Dini
distortion property:
\begin{equation}~\label{dinidistortion}
e^{-D\widetilde{\omega}(|x-y|)}\leq
\frac{(F^{-n})'(x)}{(F^{-n})'(y)}\leq e^{D\widetilde{\omega}(|x-y|)},
\quad \forall x, y\in {\mathbb R}, \;\; \forall n>0.
\end{equation}
Furthermore, let
$$
\varepsilon (t) = \left\{ \begin{array}{ll}
e^{D\widetilde{\omega}(t)}-1, & 0<t\leq 1\\
e^{D \widetilde{\omega}(1)}-1, & t>1.
\end{array}\right.
$$
Then $\varepsilon (t)>0$ is a bounded function such that
$\varepsilon (t) \to 0$ as $t\to 0^{+}$ and such that
$$
\frac{1}{1+\varepsilon (t)} \leq
\frac{(F^{-n})'(\xi)}{(F^{-n})'(\eta)} =
\frac{|F^{-n}(x+t)-F^{-n}(x)|}{|F^{-n}(x)-F^{-n}(x-t)|}\leq
1+\varepsilon (t), \quad \forall x\in {\mathbb R}, \; \forall t>0,
$$
where $\xi$ and $\eta$ are two numbers in $[0,1]$. Thus $f$ is in $\mathcal{US}$. 
 
Now we prove $f$ is uniformly absolutely continuous. Consider any homeomorphism $f^{m}: J\to I$.  Suppose $\{I_{i}\}$ are countably many pairwise disjoint open intervals in $I$. Then we have corresponding countably pairwise disjoint open intervals $\{J_{i}\}$ in $J$ such that $f^{m}(J_{i})=I_{i}$. Then we have
$$
\frac{|I_{i}|}{|I|} =\frac{|f^{m}(J_{i})|}{|f^{m}(J)|} = \frac{|(f^{m})'(\eta)||J_{i}|}{|(f^{m})'(\xi)||J|} \geq K^{-1} \frac{|J_{i}|}{|J|} 
$$
for some constant $K>0$ due to (\ref{dinidistortion}). Thus
$$
\frac{|f^{-m}(\cup_{i}I_{i}) \cap J|}{|J|}=\frac{\sum_{i} |J_{i}|}{|J|}\leq K \frac{\sum_{i} |I_{i}|}{|I|}. 
$$
This implies that $f$ is uniformly absolutely continuous. 
\end{proof}

Let ${\mathcal E}^{1+Dini}$ be the space of all $C^{1+Dini}$ expanding circle
endomorphisms of degree $d$. Then 
$$
{\mathcal E}^{1+H} \subset {\mathcal E}^{1+Dini} \subset \mathcal{UAC} \subset \mathcal{US}.
$$
We can similarly define the Teichm\"uller space $\mathcal{TE}^{1+Dini}$ as in Section~\ref{geomea}. Then we have
$$
\mathcal{TE}^{1+H} \subset \mathcal{TE}^{1+Dini} \subset \mathcal{TUAC} \subset \mathcal{TUS}.
$$

\medskip
\begin{starthmg}~\label{mainerg}
Suppose $g\in \mathcal{G}$ and $\mu$ is the unique geometric Gibbs measure associated with $g$. 
Suppose $f_{\mu}$ is uniformly absolutely continuous, that is $\tau_{\mu}=[f_{\mu}]\in \mathcal{TUAC}$. Then $f_{\mu}$ is ergodic with respect to the Lebesgue measure $|\cdot|$ on $\T$.  
\end{starthmg}

\begin{proof}
We use the same notation in the proof of Theorem A. Note that $f_{\mu}$ preserves the Lebesgue measure $|\cdot|$ on $\T$. To show that $f_{\mu}$ is ergodic, we need to show that for any Borel set $ A \ subset \T$ with $f^{-1}(A)=A$, $|A|=0$ or $ | A | = 1$.  

Suppose $A$ is a Borel set of $\T$ satisfying $f^{-1}(A)=A$ and $|A|>0$. From the Lebesgue density theorem, we have a density point $x\in A$, that is, 
$$
\lim_{t\to 0} \frac{|A\cap [x-t, x+t]|}{|[x-t,x+t]|} =1.
$$
For any $\epsilon>0$, we have a $t>0$ such that 
 $$
 \frac{|A^{c}\cap [x-t, x+t]|}{|[x-t,x+t]|} < \epsilon
 $$
 where $A^{c}=\T\setminus A$ means the complement of $A$ in $\T$. 
 
 Since $f_{\mu}$ is uniformly symmetric, we have the bounded nearby geometry for $f_{\mu}$ (see~\cite[Section 5]{JGibbs}), that is, there is a constant $D_{0}>0$ such that for any intervals $I\not=I'\in \eta_{n}$ and $I\cap I'\not= \emptyset$, $n\geq 1$, 
 $$
 \frac{1}{D_{0}} \leq \frac{|I|}{|I'|} \leq D_{0}. 
 $$
 For any fixed small $t>0$, we can find the smallest positive integer $n$ 
 such that there is a partition interval $I_{v_{n}}\in \eta_{n}$ with $I_{v_{n}} \subseteq [x-t, x+t]$. 
 Note that $I_{\sigma (v_{n})} \supset I_{v_{n}}$. Let $I_{v_{n-1}'} \in \eta_{n-1}$ be the interval having a common endpoint with $I_{\sigma (v_{n})}$. Then we have that 
 $$
 I_{v_{n}} \subseteq [x-t, x+t] \subseteq  I_{v_{n-1}'}\cup I_{\sigma (v_{n})}.
 $$
 The bounded nearby geometry implies that there is a constant $D>0$ such that 
 $$
 \frac{ |I_{v_{n-1}'}\cup I_{\sigma (v_{n})}|}{|I_{v_{n}}|} \leq D.
 $$  
 This implies that
 $$
  \frac{ |A^{c}\cap I_{v_{n}}|}{|I_{v_{n}}|} \leq    \frac{ |A^{c}\cap [x-t, x+t]|}{|I_{v_{n}}|}   \leq  D \frac{ |A^{c}\cap [x-t, x+t]|}{ |I_{v_{n-1}'}\cup I_{\sigma^{*} (v_{n})}|} 
  $$
  $$
   \leq D \frac{ |A^{c}\cap [x-t, x+t]|}{|[x-t, x+t]|}.
 $$  
 Thus we have a subsequence of integers $\{ n_{i}\}_{i=1}^{\infty}$ and a sequence of partition intervals $I_{v_{n_{i}}}\in \eta_{n_{i}}$ such that 
 $$
 \lim_{t\to \infty}   \frac{ |A\cap I_{v_{n_{i}}}|}{|I_{v_{n_{i}}}|}=1.
$$ 
 
For any fixed $n_{i}$, 
$$ 
\bigcup_{k=0}^{\infty} f^{-kn_{i}} (I_{v_{n_{i}}}) 
= I_{v_{n_{i}}} \bigcup \bigcup_{k=1}^{\infty} \bigcup_{v_{n_{i}}^{k}\not=v_{n_{i}}, \cdots, v_{n_{i}}^{1}\not= v_{n_{i}}} I_{v_{n_{i}}^{k}\cdots v_{n_{i}}^{1}v_{n_{i}}}
$$
$$
=\T\setminus \Big( \bigcup_{k=1}^{\infty} \bigcup_{v_{n_{i}}^{k}\not= v_{n_{i}}, \cdots, v_{n_{i}}^{1} \not=v_{n_{i}}} I_{v_{n_{i}}^{k}\cdots v_{n_{i}}^{1}}\Big). 
$$
As we have proved before (see (\ref{zero1})),  
$$
\mu_{0} \Big(\bigcup_{k=0}^{\infty} f^{-kn_{i}} (I_{v_{n_{i}}})\Big) =\lim_{k\to \infty} \Big(
|I_{v_{n_{i}}}| +\sum_{l=1}^{k} \sum_{v_{n_{i}}^{l}\not=v_{n_{i}}, \cdots, v_{n_{i}}^{1}\not= v_{n_{i}}} |I_{v_{n_{i}}^{l}\cdots v_{n_{i}}^{1}v_{n_{i}}}|\Big)
$$
$$
=1- \lim_{k\to \infty}  \sum_{l=1}^{k} \sum_{v_{n_{i}}^{l}\not= v_{n_{i}}, \cdots, v_{n_{i}}^{1} \not=v_{n_{i}}} I_{v_{n_{i}}^{l}\cdots v_{n_{i}}^{1}} =1-0=1.
$$

For any $\epsilon >0$, since $f_{\mu}$ is uniformly absolutely continuous, we have a $\delta>0$ 
such that for any $n_{i}$ and any countably many pairwise disjoint open intervals $\{J_{j}\}$ in $I_{v_{n_{i}}}$ with 
$$
\frac{ \sum_{j} |J_{j}|}{|I_{v_{n_{i}}}|}<\delta, 
$$
$$
\frac{|f^{-kn_{i}} (\cup_{j}J_{j})\cap I_{v_{n_{i}}^{k}\cdots v_{n_{i}}^{1}v_{n_{i}}}|}{|I_{v_{n_{i}}^{k}\cdots v_{n_{i}}^{1}v_{n_{i}}}|} <\epsilon.
$$
Since 
$$
 \lim_{t\to \infty}   \frac{ |A^{c}\cap I_{v_{n_{i}}}|}{|I_{v_{n_{i}}}|}=0,
$$ 
we can find a $i_{0}>0$ such that 
$$
\frac{ |A^{c}\cap I_{v_{n_{i}}}|}{|I_{v_{n_{i}}}|}<\delta
$$
for all $i>i_{0}$. For each $i>i_{0}$, we have an open set $B \supset A^{c}$. 
$$
\frac{ |B\cap I_{v_{n_{i}}}|}{|I_{v_{n_{i}}}|}<\delta.
$$
The interior of $B\cap I_{v_{n_{i}}}$ is the union of countably many pairwise disjoint intervals in $I_{v_{n_{i}}}$. Thus we have that
$$
\frac{|f^{-kn_{i}} (A^{c})\cap I_{v_{n_{i}}^{k}\cdots v_{n_{i}}^{1}v_{n_{i}}}|}{|I_{v_{n_{i}}^{k}\cdots v_{n_{i}}^{1}v_{n_{i}}}|} 
\leq 
\frac{|f^{-kn_{i}} (B)\cap I_{v_{n_{i}}^{k}\cdots v_{n_{i}}^{1}v_{n_{i}}}|}{|I_{v_{n_{i}}^{k}\cdots v_{n_{i}}^{1}v_{n_{i}}}|} <\epsilon.
$$
This implies that 
$$
\frac{|f^{-kn_{i}} (A^{c})\cap \cup_{l=0}^{k} f^{-ln_{i}} (I_{v_{n_{i}}}) |}{|\cup_{l=0}^{k} f^{-ln_{i}} (I_{v_{n_{i}}}) |} <\epsilon.
$$
Let $k\to \infty$, we have that
$$
|A^{c}| =|A^{c}\cap \cup_{k=0}^{\infty} f^{-kn_{i}} (I_{v_{n_{i}}})| <\epsilon
$$
since $f^{-kn_{i}}(A^{c})=A^{c}$ for all $k$ and $|\cup_{k=0}^{\infty} f^{-kn_{i}} (I_{v_{n_{i}}})|=1$.
This implies that $|A^{c}|=0$, that is, $|A|=1$.  We proved the theorem.
 
\end{proof}

\medskip
\begin{corollary}
Suppose $g\in \mathcal{G}$ and $\mu$ is the unique geometric Gibbs measure associated with $g$. 
Suppose $f_{\mu}$ is uniformly absolutely continuous. Then the geometric Gibbs measure $\mu$ is ergodic. 
\end{corollary}

\medskip
\begin{remark}~\label{uwp}
In general, a uniformly symmetric circle endomorphism $f\in\mathcal{US}$ is not absolutely continuous and hence not uniformly absolutely continuous. For each $f\in\mathcal{US}$, let $F$ be its lift to the real line $\R$. Suppose $\widehat{F}$ is a quasiconformal extension of $F$ to the upper-half plane $\UH$. This yields a sequence of Beltrami coefficients \(
\{b_{\widehat{F}^{-n}}\}_{n=0}^{\infty}\).
By Theorem~\ref{ccusthm}, $f$ is uniformly symmetric if and only if there exists an extension $\widehat{F}$ such that $b_{\widehat{F}^{-n}}(z)\to 0$ as $z$ approaches the boundary $\R$, uniformly in $n$.
The notion of a uniformly absolutely continuous circle endomorphism (Definition~\ref{uac}) relates to the Weil–Petersson (WP) property; see~\cite{Bi,Sh}. We define the \emph{uniformly WP} condition for $f$ by requiring
\begin{equation}~\label{uwpc}
\sup_{n\ge 1}\;\int_{\UH} \frac{|b_{\widehat{F}^{-n}}(z)|^{2}}{y^{2}} dx dy< \infty,\quad z=x+iy\in \UH.
\end{equation}
According to (\ref{ceq}), (\ref{uwpc}) is equivalent to
\begin{equation}~\label{uwpcc}
\sup_{n\ge 1}\;\int_{\UH} \frac{|b_{\widetilde{H}} (d^{-n}z) -b_{\widetilde{H}}(z)|}{y^{2}} dx dy<\infty, \quad z=x+iy\in \UH.
\end{equation}
In general, $H=\widetilde{H}|{\mathbb R}$ is not absolutely continuous. We still consider the following problem. 

\medskip
\begin{problem}~\label{prob3}
Is a uniformly WP circle endomorphism $f$ uniformly absolutely continuous?
\end{problem}

\end{remark}



\bibliographystyle{amsalpha}

\end{document}